\documentclass[10pt]{amsart}
\pdfoutput=1
\usepackage{eufrak}
\usepackage{amsfonts}
\usepackage{float}
\usepackage{amsmath,amssymb}
\usepackage{lscape}
\usepackage[pdftex]{graphicx, color}
\usepackage[all]{xy}

\def\qed{\hfill $\Box$}
\def\proof{\noindent {\sl Proof} :\;  }

\newcommand{\bP}{\mathbb{P}}
\newcommand{\C}{\mathbb{C}}
\newcommand{\Z}{\mathbb{Z}}

\newcommand{\R}{\mathbb{R}}

\def\qed{\hfill $\Box$}
\def\proof{\noindent {\sl Proof} :\;  }

\def\rd{\partial}

\def\bx{\mbox{\boldmath $x$}}

\def\bu{\mbox{\boldmath $u$}}
\def\bv{\mbox{\boldmath $v$}}

\def\ba{\mbox{\boldmath $a$}}
\def\bb{\mbox{\boldmath $b$}}
\def\bc{\mbox{\boldmath $c$}}

\def\be{\mbox{\boldmath $e$}}

\def\bm{\mbox{\boldmath $m$}}

\def\bp{\mbox{\boldmath $p$}}

\def\b0{\mbox{\boldmath $0$}}

\def\bbx{\mbox{\tiny$\bx$}}

\def\bba{\mbox{\tiny$\ba$}}
\def\bbp{\mbox{\tiny$\bp$}}

\newtheorem{thm}{\bf Theorem}[section]

\newtheorem{lem}[thm]{\bf Lemma}
\newtheorem{prop}[thm]{\bf Proposition}
\newtheorem{dfn}[thm]{\bf Definition}

\newtheorem{rem}[thm]{\bf Remark}
\newtheorem{exam}[thm]{\bf Example}

\begin{document}
\title[The dually flat structure for singular models]
{The dually flat structure for singular models}
\author[N.~Nakajima]{Naomichi Nakajima}
\address[N.~Nakajima]{M2, Graduate School of Information Science and Technology, Hokkaido University,
Sapporo 060-0814, Japan}
\email{nakajima-n@ist.hokudai.ac.jp}
\author[T.~Ohmoto]{Toru Ohmoto}
\address[T.~Ohmoto]{Department of Mathematics,
Faculty of Science,  Global Station of Bigdata and cybersecurity (GiCORE-GSB), Hokkaido University,
Sapporo 060-0810, Japan}
\email{ohmoto@math.sci.hokudai.ac.jp}
%
%
\dedicatory{}

%
\begin{abstract}
The {\it dually flat structure} introduced by Amari-Nagaoka is highlighted in information geometry and related fields. In practical applications,  however,  the underlying pseudo-Riemannian metric may often be degenerate, 
and such an excellent geometric structure is rarely defined on the entire space. 
To fix this trouble, in the present paper, we propose a novel generalization of the dually flat structure for a certain class of singular models from the viewpoint of {\em Lagrange and Legendre singularity theory} --  we introduce a {\em quasi-Hessian manifold} endowed with a possibly degenerate metric and a particular symmetric cubic tensor, which exceeds the concept of statistical manifolds and is adapted to the theory of (weak) contrast functions.  In particular, we establish Amari-Nagaoka's extended Pythagorean theorem and projection theorem in this general setup, and consequently, most of applications of these theorems are suitably justified even for such singular cases. This work is motivated by various interests with different backgrounds from Frobenius structure in mathematical physics to Deep Learning in data science. 

\end{abstract}

\maketitle

{\small 
\setcounter{tocdepth}{2}
\tableofcontents
}

\newpage

\section{Introduction}
\label{intro}
The {\it dually flat structure} is highlighted in information geometry  
-- it brings a united geometric insight on various fields such as statistical science, convex optimizations, (quantum) information theory, and so on (Amari-Nagaoka \cite{AmariNagaoka00}, Amari \cite{Amari16}, Chentsov \cite{Chentsov}). 
This is also essentially the same as the {\em Hessian structure} in affine differential geometry (Shima \cite{Shima}). 
On a $C^\infty$-manifold $M$, 
a dually flat structure is a triplet $(h, \nabla, \nabla^*)$ 
where $h$ is a pseudo-Riemannian metric (i.e., non-degenerate symmetric $(0,2)$-tensor) 
and $\nabla$ and $\nabla^*$ are flat affine connections on $M$ satisfying certain properties; the most particular feature is that the metric is locally given by the Hessian matrix of some potential function in $\nabla$-affine coordinates. In practical applications, however, the Hessian matrix may often be degenerate along some locus $\Sigma$ of points in $M$, and then, strictly speaking, the differential geometric method can not be directly applied. 
We call such a space a {\em singular model}, roughly. 
In the present paper, we propose a novel generalization of the dually flat structure for a certain class of singular models from the viewpoint of contact geometry and singularity theory. 
This provides a new framework for general hierarchical structures -- we introduce a {\em quasi-Hessian manifold} $M$ endowed with a possibly degenerate quadratic tensor $h$ and a particular symmetric cubic tensor $C$, that exceeds the concept of statistical manifolds and very fits with the theory of contrast functions due to Eguchi \cite{Eguchi}.  
In fact, such $M$ naturally possesses a {\em canonical divergence} $\mathcal{D}: M\times M \to \R$, 
which is a {\em weak} contrast function compatible with $h$ and $C$ (Theorem \ref{contrast}). 
The key is the {\em Legendre duality}, which {\em does} exist even under the presence of the degeneracy locus $\Sigma$ of $h$.  In spite of no metric $h$ and no connection $\nabla$ available (!), we generalize in a natural way the Amari-Chentsov cubic tensor $\nabla h$ to a symmetric tensor $C$ defined on the entire space $M$ (that is possible even in case that $M=\Sigma$), and especially we establish Amari-Nagaoka's {\em extended Pythagorean theorem} and {\em projection theorem} in this setup (Theorems \ref{pytha_thm}, \ref{proj_thm}). Consequently, in principle, most of applications of these theorems are suitably justified even for such degenerate cases. 

As the first observation, we see that if the Hessian of a potential is degenerate, the graph of the dual potential (i.e., the Legendre transform of the potential) is no longer a submanifold but  a {\em wavefront} having singularities branched along its {\em caustics}. 
More generally, our quasi-Hessian manifold $M$ is locally accompanied with two kinds of wavefronts, later called the {\em $e/m$-wavefronts}, as an alternate to the pair of a convex potential and its dual. 
To grasp the point, it would be helpful to refer to Fig.1 and Fig.2 in \S \ref{sec:3-1} in advance. 
Those two wavefronts are mutually tied by the Legendre duality in a strict sense, and also they have  `height functions' (i.e., projections to the $z$ and $z'$-axes, respectively) by using which we generalize the Bregman divergence to our divergence $\mathcal{D}$. 
Further, we associate the pair of {\it coherent tangent bundles} 
$$(E, \nabla^E, \Phi: TM \to E) \;\; \mbox{and}\;\;  (E', \nabla^{E'}, \Phi': TM \to E'),$$
where each of $E$ and $E'$ is a vector bundle on $M$  that is an alternative to the tangent bundle $TM$ and equipped with a flat connection and a bundle map from $TM$ measuring the degeneration of $h$, and $E$ and $E'$ are dual to each other. 
Intuitively, the coherent tangent bundles come from the splitting of the standard contact structure into the directions of the $e/m$-wavefront projections (see \S \ref{sec:2-1}, \S \ref{sec:3-1} and \S \ref{sec:3-2}). 
The role of $\nabla$ and $\nabla^*$ on a dually flat space is undertaken by $\nabla^E$ and $\nabla^{E'}$ of those vector bundles on our singular model $M$, and then the new cubic tensor $C$ of $TM$ is defined by using these connections  through $\Phi$ and $\Phi'$ (\S \ref{sec:3-4}). 
Originally,  the notion of coherent tangent bundles has been introduced for studying Riemannian geometry of singular wavefronts by Saji-Umehara-Yamada \cite{SUY}, and here we borrow an affine geometry version. 
 In the case that $h$ is non-degenerate ($E=TM$, $E'=T^*M$),  the dually flat structure in the original form is naturally recovered. 
 
Singularities of caustics and wavefronts have thoroughly been investigated in  {\it Lagrange and Legendre singularity theory} (initiated by Arnol'd, Zakalyukin and H\"ormander) in relation with a broad range of subjects such as classical mechanics, thermodynamics, geometric optics, Fourier integral operators, control theory, catastrophe theory and so on \cite{Arnold89, AGV, Arnold90, Arnold92, Izumiya_Ishikawa, ICRT, PS}. We bring several techniques or concepts in this theory into information geometry, that may suggest new directions in both theory and application. 
In fact, the present paper is motivated by various interests from different backgrounds:
\begin{itemize}
\item[-]  
A typical example of quasi-Hessian manifolds is a general affine hypersurface $M$ in $\R^{n+1}$. 
It possesses {\it mixed geometry} with changing metric-type -- that goes back to Darboux and others  
dealing with a rich geometry of {\em parabolic curve $\Sigma$} (the curve of inflection points) separating elliptic and hyperbolic domains on a surface in $\R^3$  \cite{Arnold92, Arnold90, SKDO}. 
That is also related to {\em nonconvex} optimization and variational problems \cite{Ekeland}. 
\item[-] 
Any Lagrange submanifold of $\R^{2n}$ is a quasi-Hessian manifold. 
If it is flat, then the cubic tensor $C$ satisfies the so-called WDVV equation, which mainly appears in topological field theory, that yields a version of Frobenius manifold-like structure \cite{Kito, Totaro, Hertling}. 
That is also related to geometry of Poisson manifolds and paraK\"ahler structure \cite{Boucetta19, Combe}.  
\item[-] 
In statistical inference, 
any curved exponential family produces a quasi-Hessian manifold, which represents the $\nabla^*$-extrinsic geometry in the ambient family. For instance, it is useful for studying catastrophe phenomena in root selections of the maximum likelihood equation \cite{Small}. Almost all statistical learning machines including deep neural networks allow degeneration of Fisher-Rao matrices \cite{Amari16, FK, Watanabe}, to which we are seeking for a new approach. 
\end{itemize}

In the present paper, we will focus only on basic ideas and write them plainly in a self-contained manner as much as possible --  most arguments are elementary and checked by direct computations, and we will NOT enter into any detail of singularity theory here. Therefore, perhaps, this paper would be readable enough for anyone with various background. 
Nevertheless, we believe that this paper contains some new observations in this field. 
The rest of this paper is organized as follows. 
In \S 2 we give a brief summary on some basics in contact geometry and the dually flat structure. 
In \S 3, after reviewing the definition of Lagrange and Legendre singularities, we introduce quasi-Hessian manifolds. 
In \S 4, the associated canonical divergence will be discussed; 
the extended Pythagorean theorem and projection theorem are presented in our setting 
and also we give a relation with contrast functions. 
In \S 5, we pick up some possible applications and open questions. 

Throughout, bold letters denote column vectors, e.g., $\bx=(x_1, \cdots, x_n)^T$, 
and the notation with prime $\bx'$ simply means to distinguish from $\bx$ 
(not mean any operation like differential or transpose). 
Also we let $\frac{\rd f}{\rd \bx}$ denote $(\frac{\rd f}{\rd x_1}, \cdots, \frac{\rd f}{\rd x_n})^T$ for short as usual. 
We assume that manifolds and maps are of class $C^\infty$, for the simplicity.

The authors are partly supported by GiCORE-GSB, Hokkaido University, 
and JSPS KAKENHI Grant Numbers JP17H06128 and JP18K18714.

\

\section{Dually Flat Structure}\label{sec:2}

\subsection{Contact geometry and Legendre duality}\label{sec:2-1} 
To begin with, we summarize a minimal set of basic knowledge in contact geometry which will be used throughout this paper. 
As best references, we recommend Chap.18-22 of Arnol'd et al \cite{AGV},  Appendix of Arnol'd \cite{Arnold89} and Izumiya-Ishikawa \cite{Izumiya_Ishikawa}. 

A {\it contact manifold} is a $(2n+1)$-dimensional manifold $N$ endowed with 
a maximally non-integrable hyperplane field $\xi_p \subset T_pN\; (p \in N)$, 
i.e.,  $\xi$ is locally expressed by the kernel of a $1$-form $\theta$ satisfying the non-degeneracy condition $\theta\wedge (d\theta)^n \not=0$. This field $\xi$ is called a {\em contact structure} on $N$. 
The most important example is the {\it standard contact space}  $\R^{2n+1}$; 
it is the {\em $1$-jet space} of functions on the affine $n$-space 
$$\R^{2n+1}=J^1(\R^n, \R)=T^*\R^n\times \R,$$
where the $1$-jet of a function $f$ at a point $\ba$ means the Taylor coefficients at $\ba$ of order $\le 1$, i.e., $(df(\ba), f(\ba)) \in T_{\bba}^*\R^n\times \R$. 
The contact structure is given by the $1$-form 
$$\theta=dz-\bp ^Td\bx=dz - \sum_{i=1}^n p_idx_i,$$ 
called the {\em standard contact form}, where 
$z$ is the last coordinate, and $\bx$ and $\bp$ denote, respectively,  
the base and the fiber coordinates of the cotangent bundle $T^*\R^n$ (we always write coordinates in this order). 
We often write $\R^n_{\bbx}$ and $\R^n_{\bbp}$ in order to distinguish them. 
Note that the standard contact structure relies on the affine structure of the base space, 
not on the choice of coordinates $\bx$.  

The famous {\em Darboux theorem} tells that {\em the contact structure is locally unique}; namely, for any contact manifold $N$, we can always find a system of local coordinates around any point $p \in N$,  in which the contact structure is presented by the standard one. 

A {\it Legendre submanifold} $L$ of a $(2n+1)$-dimensional contact manifold $(N, \xi)$ is 
an $n$-dimensional integral manifold of the field $\xi$, i.e.,  
$T_pL \subset \xi_p$ for every $p \in L$. 
It is easy to see that in the standard $\R^{2n+1}$, 
the graph of $(df, f)$ of a function $z=f(\bx)$ 
$$L=\left\{\,(\bx, \bp, z) \in \R^{2n+1}\, \bigg|\, \bp= \frac{\rd f}{\rd \bx}, \; z=f(\bx)\, \right\}
\eqno{(1)}$$
is a Legendre submanifold, 
and conversely, every Legendre submanifold 
which is diffeomorphically projected to $\R^n_{\bbx}$ is expressed in this form (1). 

The {\em standard symplectic structure} of $T^*\R^n$ is defined by the non-degenerate closed $2$-form 
$$\omega=\sum_{i=1}^n dx_i\wedge dp_i.$$ 
A {\em Lagrange submanifold} is an $n$-dimensional submanifold over which $\omega$ vanishes. 
A typical example is the graph of $df$, 
i.e., the image of $L$ of the form (1) via projection along the $z$-axis. 
Any Lagrange submanifold of $T^*\R^n$ is always {\em locally} liftable to a Legendre submanifold of $T^*\R^n\times \R$ uniquely up to a transition parallel to the $z$-axis. If it is entirely liftable, 
then we call it an {\em exact} Lagrange submanifold.

A {\it Legendre fibration} $\pi: N \to B$ is a fiber bundle whose total space is a $(2n+1)$-dimensional contact manifold $N$, 
the base is an $(n+1)$-dimensional manifold $B$ and every fiber $\pi^{-1}(x)$ is Legendrian. 
The most typical example is the projection from the standard space 
$$\pi: \R^{2n+1} \to \R^n_{\bbx}\times \R, \;\; (\bx, \bp, z) \mapsto (\bx, z).$$ 
Every Legendre fibration is locally described in this typical form with suitable local coordinates.  

The {\em Legendre duality} is described as follows. 
Consider the transformation $\mathcal{L}: \R^{2n+1} \to \R^{2n+1}$ given by 
$$(\bx', \bp', z')=\mathcal{L}(\bx, \bp, z) := (\bp, \bx, \bp^T\bx-z).$$
It is a {\it contactomorphism}, i.e., 
$\mathcal{L}$ preserves the contact hyperplane fields $\xi$; 
indeed,  $\mathcal{L}^*\theta=\mathcal{L}^*(dz'-\bp'^Td\bx')=- \theta$. 
Put 
$$\pi':=\pi\circ \mathcal{L}: \R^{2n+1} \to \R^n_{\bbp}\times \R, \; (\bx, \bp, z) \mapsto (\bp, \bp^T\bx-z),$$ 
which is also a Legendre fibration. 
Then, the {\it double fibration structure} of the standard contact space is defined as the following diagram: 
$$\xymatrix{\R^{n}_{\bbx}\times \R & \R^{2n+1} \ar[l]_{\;\;\pi} \ar[r]^{\pi'\;\;} & \R^{n}_{\bbp}\times \R} \eqno{(dL)}$$
Let $\Pi: \R^n\times \R \to \R^n$ be the projection to the first factor 
and put 
$$\pi_1=\Pi\circ \pi: \R^{2n+1} \to \R^{n}_{\bbx}, \quad \pi'_1=\Pi\circ \pi': \R^{2n+1} \to \R^{n}_{\bbp}.$$ 
Let $L$ be a Legendre submanifold of $\R^{2n+1}$. 
In this paper, we call $L$ a {\em regular model} if $L$ is diffeomorphic to 
some open subsets $U\subset \R^n_{\bbx}$ and $V\subset \R^n_{\bbp}$ via projections $\pi_1$ and $\pi'_1$, respectively.  
Equivalently, 
there exsit functions $z=f(\bx)$ on $U$ and $z'=\varphi(\bp)$ on $V$ such that 
\begin{itemize}
\item[-\;\;] $L  \subset \R^{2n+1}=J^1(\R^n_{\bbx},\R)$ is the graph of $(df, f)$;  
\item[-\;\;] $\mathcal{L}(L)\subset \R^{2n+1}=J^1(\R^n_{\bbp},\R)$ is the graph of $(d\varphi, \varphi)$.
\end{itemize}
Then we have 
$$\bp=\frac{\rd f}{\rd \bx}, \;\; \bx=\frac{\rd \varphi}{\rd \bp}, \quad f(\bx)+\varphi(\bp) - \bp^T\bx=0.$$
The coordinate change is the gradient map $U \to V$, $\bx \mapsto \bp=\frac{\rd f}{\rd \bbx}$. 
It is diffeomorphic, thus the Hessian matrix of $f(\bx)$ is non-degenerate at every $\bx \in U$.  
Here, the inverse map $V \to U$ is given by $\bp \mapsto \bx=\frac{\rd \varphi}{\rd \bbp}$,  
and its Hessian matrix is the inverse of that of $f(\bx)$. 
We say that $z'=\varphi(\bp)$ is the {\em Legendre transform} of $z=f(\bx)$ and vice-versa. 
We call $f$ a {\em potential function} and $\varphi$ its {\em dual potential}. 
This correspondence is the Legendre duality. It is very common in, e.g., convex analysis: 
if $z=f(\bx)$ is strictly convex, then $z'=\varphi(\bp)$ is also (see Remark \ref{proj}). 

An {\em affine Legendre equivalence}, a new terminology introduced in this paper, 
is defined by an affine transformation 
$\mathcal{L}_{F}: \R^{2n+1} \to \R^{2n+1}$ of the form 
$$\mathcal{L}_{F} (\bx, \bp, z) = (A\bx+ \bb, \;  A'\bp+ \bb', \; z+\bc^T\bx+ d)$$
together with affine transformations 
$$F(\bx, z)=(A\bx+ \bb,  z+\bc^T\bx+ d),$$ 
$$F^*(\bp, z')= (A'\bp+ \bb', z'+\bc'^{T}\bp+ d'),$$ 
where $A$ is invertible and 
$$A'=(A^T)^{-1}, \;\; \bb'=A'\bc, \;\; \bb=A\bc', \;\; d'=\bb'^T\bb-d.$$ 
Note that  $F$ (or $F^*$) determines $\mathcal{L}_F$. 
It is easy to see that $\mathcal{L}_{F}$ preserves the contact form and the double fibrations $(dL)$, i.e., 
$\mathcal{L}_{F}^*\theta=\theta$ and the following diagram commutes: 
$$
\xymatrix{
\R^{n}_{\bbx}\times \R_z \ar[d]_F^{\simeq}& \R^{2n+1} \ar[l]_{\;\;\pi} \ar[r]^{\pi'\;\;}\ar[d]^{\mathcal{L}_F}_{\simeq} & \R^{n}_{\bbp}\times \R_{z'} \ar[d]^{F^*}_{\simeq}\\
\R^{n}_{\bbx}\times \R_z & \R^{2n+1} \ar[l]_{\;\;\pi} \ar[r]^{\pi'\;\;} & \R^{n}_{\bbp}\times \R_{z'}
}
$$

\begin{dfn}\upshape
We say that two Legendre submanifolds $L_1, L_2$ of $\R^{2n+1}$ are {\em affine Legendre equivalent}  if 
there exists some $\mathcal{L}_{F}$ which identifies $L_1$ with $L_2$. 
\end{dfn}

\begin{rem}\label{proj}
\upshape
{\bf (Projective duality)} 
The Legendre duality is an {\em affine expression} of the projective duality. 
We denote by $\bP^{n+1}\, (:=\R P^{n+1})$ the real projective space of dimension $n+1$ 
and by $\bP^{n+1*}\, (:=\R P^{n+1*})$ the dual projective space. 
Let $N$ denote the incidence submanifold of $\bP^{n+1}\times \bP^{n+1*}$ which consists of pairs $(p, \lambda)$ with $p \in \lambda$, i.e., $N$ is a codimension one submanifold ($\dim N=2n+1$) defined by 
$$p_0x_0+p_1x_1+\cdots +p_{n+1}x_{n+1}=0$$
for $p=[x_0:\cdots: x_{n+1}] \in \bP^{n+1}$ and  $\lambda =[p_0:\cdots: p_{n+1}]\in \bP^{n+1*}$.  
Note that $N$ is naturally identified with the projective cotangent bundle $PT^*\bP^{n+1}\, (=PT^*\bP^{n+1*})$, and thus $N$ becomes a contact manifold  \cite[\S 20.1]{AGV}. 
Consider the open subset $O_N$ of $N$ defined by $x_{n+1}\not=0$ and $p_{0}\not=0$. 
We may set $x_{n+1}=p_{0}=-1$, and put  $z=x_0$ and $z'=p_{n+1}$, then the above equation is rewritten as 
$$z+z'-\bp^T\bx=0.$$  
Clearly, $O_N$ has two systems of coordinates, $(\bx, \bp, z)$ and $(\bp, \bx, z')$, 
and the coordinate change between them is just the above $\mathcal{L}: \R^{2n+1} \to \R^{2n+1}$ 
preserving the contact structure of $N$. 
In projective geometry, 
the double Legendre fibration 
$$\xymatrix{
\bP^{n+1} & N \ar[l]_{\;\;\;\pi} \ar[r]^{\pi'\;\;\;} & \bP^{n+1*}
} $$
expresses the {\em duality principle on points and hyperplanes}, 
where $\pi$ and $\pi'$ are projections of the projective cotangent bundles. 
Restrict this diagram to $O_N$ 
and identify $O_N$ with $\R^{2n+1}$ using coordinates $(\bx, \bp, z)$, 
we get the diagram $(dL)$. 
For instance, in case of $n=1$, consider a parameterized plane curve 
$$C: (x, z):=(x, f(x))\in \R^2 \subset \bP^2.$$ 
Then its projective dual is the following curve consisting of the tangent lines: 
$$\textstyle C^*: (p, z')=(\frac{df}{dx}(x), x\frac{df}{dx}(x)-f(x)) \in \R^{2*} \subset \bP^{2*}.$$ 
If $C$ is convex, then $C^*$ is also.  
If $C$ has an inflection point, e.g., $f(x)=\frac{1}{3}x^3+\cdots$, 
then $C^*$ has a cusp at the corresponding point, $(p, z')=(x^2+\cdots, \frac{2}{3}x^3+\cdots)$, 
and therefore, 
$C^*$ is locally the graph of a bi-valued function, $z'=\pm \frac{2}{3}p^{3/2}+\cdots\, (p\ge 0)$. 
\end{rem}

\subsection{Dually flat structure}\label{sec:2-2}
Let $L$ be a Legendre submanifold of a regular model with potential function $z=f(\bx)$. 
The Hessian matrix 
$$H(p)=\left[\frac{\rd^2 f}{\rd x_i\rd x_j}(\pi_1(p))\right] \quad (p \in L)$$ 
is invertible, thus it defines a pseudo-Riemannian metric $h$ on $L$, called the {\em Hessian metric associated to $f$}. 
Additionally, through the projections $\pi_1$ and $\pi'_1$, 
the fixed affine structures of $\R^n_{\bbx}$ and $\R^n_{\bbp}$ induce two different flat affine connections $\nabla, \nabla^*$ on $L$,  respectively. 

\begin{dfn}\label{duallyflat}
\upshape 
{(\cite{AmariNagaoka00})} 
The triplet $(h, \nabla, \nabla^*)$ is called the {\em dually flat structure} on a regular model $L$. 
\end{dfn}

Note that $(h, \nabla, \nabla^*)$ is preserved under affine Legendre equivalence; indeed, $\mathcal{L}_F$ induces affine transformations of $\R^n_{\bbx}$ and $\R^n_{\bbp}$ and simply adds a linear function $\bc^T\bx+d$ to the potential $z=f(\bx)$. 

The dually flat structure is traditionally introduced in terms of differential geometry in an intrinsic way. 
We briefly summarize it below, see \cite{AmariNagaoka00, Amari16, Matsuzoe, Shima} for the detail. 

A {\it statistical manifold} is a pseudo-Riemannian manifold $(M, h)$ equipped with  
a torsion-free affine connection $\nabla$ being compatible with $h$, i.e., 
the cubic tensor $T:=\nabla h$ is totally symmetric: 
$$(\nabla_Xh)(Y, Z)=(\nabla_Yh)(X, Z)$$ 
for vector fields $X, Y$ and $Z$. 
Equivalently \cite[p.306]{Matsuzoe},  
a stastistical manifold may also be defined as a manifold endowed  
with a pseudo-Riemannian metric $h$ and a totally symmetric $(0,3)$-tensor $T$ (due to Lauritzen), that is also described within the theory of contrast functions in \cite{Eguchi}. 
The {\it dual connection $\nabla^*$} (with respect to $h$) is defined by 
$$Xh(Y, Z)=h(\nabla_XY, Z)+h(Y, \nabla_X^*Z),$$ 
and then $\nabla^*$ is torsion-free and $\nabla^*h$ is also symmetric. 
Furthermore, if $\nabla$ is flat (i.e., torsion-free and curvature-free),  then $\nabla^*$ is also. 
Such a statistical manifold with flat connections is called a {\it dually flat manifold} \cite{AmariNagaoka00, Amari16} or 
a {\em Hessian manifold}  \cite{Shima, Matsuzoe}. 
The most notable characteristic of a dually flat manifold is that locally it holds that  
$$h=\nabla df$$ 
for some local potential $f$. 
In other words, the metric $h$ is expressed by  
the non-degenerate Hessian matrix of $z=f(\bx)$ in $\nabla$-affine local coordinates $\bx$. 
Moreover, the $\nabla^*$-affine coordinates $\bp$ are then given by $\bp=\frac{\rd f}{\rd \bbx}(\bx)$. 

This local expression of a dually flat manifold $M$ exactly provides a regular model $L$ in $\R^{2n+1}$, 
the graph of $1$-jet of a local potential,  
equipped with the dually flat structure in the sense of Definition \ref{duallyflat}. 
Such a regular model $L$ is {\em uniquely determined up to affine Legendre equivalence}. 
To see this precisely, suppose that the metric $h$ is locally expressed by the Hessian matrices $H_\alpha$ and $H_\beta$ of two potential functions $f_\alpha(\bx^\alpha)$ and $f_\beta(\bx^\beta)$   
in different $\nabla$-affine local coordinates, respectively. 
Here, let $(U_\alpha, \bx^\alpha=(x^\alpha_1, \cdots, x^\alpha_n)^T)$ and 
$(U_\beta, \bx^\beta=(x^\beta_1, \cdots, x^\beta_n)^T)$ 
denote the charts with $U_\alpha, U_\beta \subset M$, $U_\alpha \cap U_\beta\not=\emptyset$. 
By definition, there is an affine transformation 
$$\bx^\beta=\psi(\bx^\alpha)=A\bx^\alpha + \bb.$$ 
By the assumption,  it holds that 
$A^T H_\beta(p) A=H_\alpha(p)$ for every $p \in U_\alpha \cap U_\beta$, 
thus any second partial derivatives of the composite function $f_\beta\circ \psi(\bx^\alpha)$ coincide with those of $f_\alpha(\bx^\alpha)$. 
Namely, these two functions are the same up to some linear term: 
$$f_\beta\circ \psi(\bx^\alpha)=f_\alpha(\bx^\alpha)+\bc^T\bx^\alpha+d.$$ 
Then the affine transformation 
$$F(\bx^\alpha, z):=(A\bx^\alpha+ \bb,  z+\bc^T\bx^\alpha+ d)$$ 
sends the graph of $z=f_\alpha(\bx^\alpha)$ to the graph of $z=f_\beta(\bx^\beta)$. 
Hence, the corresponding affine Legendre equivalence $\mathcal{L}_F$ identifies 
two regular models,  $L_\alpha \subset \R^{2n+1}$ defined by $(df_\alpha, f_\alpha)$ 
and $L_\beta \subset \R^{2n+1}$ defined by $(df_\beta, f_\beta)$, 
on the overlap. 

Actually, less noticed, though, this simple observation says that any dually flat manifold is an affine manifold having   
an atlas $\{(U_\alpha, \bx^\alpha)\}_{\alpha \in \Lambda}$ with affine coordinate changes $\psi_{\alpha}^\beta\, (\alpha, \beta \in \Lambda)$ so that it is additionally equipped with local potentials $\{f_\alpha\}_{\alpha \in \Lambda}$ whose graphs are glued by affine transformations $F_\alpha^\beta$ of the above form. 
The affine structure gives the flat connection $\nabla$, and local potentials restore the Hessian metric $h$ by gluing $\{H_\alpha\}_{\alpha \in \Lambda}$. At the level of Legendre submanifolds given by $1$-jets of local potentials, notice again that $\mathcal{L}_F$ preserves $(h, \nabla, \nabla^*)$. 
Therefore, we may rephrase the above statement in the following way: 

\begin{prop}\label{Hessian_manifold}
Any dually flat or Hessian manifold is an affine manifold made up by gluing several regular models $L_\alpha$ in $\R^{2n+1}$ via affine Legendre equivalence. 
The metric $h$ and the pair of affine connections $\nabla$ and $\nabla^*$ are reconstructed by the dually flat structures of $L_\alpha$  in the sense of Definition \ref{duallyflat}. 
\end{prop}

This gluing construction will be generalized later to introduce our quasi-Hessian manifolds (Definition \ref{quasi-Hessian_manifold} in \S \ref{sec:3-3}). 

\begin{rem}
\upshape
Since each gluing map acts also on a neighborhood of a regular model in $\R^{2n+1}$, 
the gluing construction yields a dually flat manifold as a Legendre submanifold of some ambient contact manifold (also it produces a Lagrange submanifold of some symplectic manifold). 
Let $(M, h, \nabla, \nabla^*)$ be a dually flat manifold, and 
suppose that there exists a global potential $f: M \to \R$ with $h=\nabla df$. 
Take local charts $U_\alpha$ of $\nabla$-flat coordinates, then local potentials $f|_{U_\alpha}$ define regular models $L_\alpha$ 
in $J^1(U_\alpha, \R)=T^*U_\alpha \times \R$, and they are glued together by affine Legendre equivalence of the form $\mathcal{L}_F$ with $\bc=0$ and $d=0$. 
Conversely, gluing local models by this special kind of affine Legendre equivalences yields 
a dually flat manifold with a global potential. 
As a weaker situation, suppose that there exists a closed $1$-form $\sigma$ with $h=\nabla \sigma$; 
then $M$ is said to be of {\em Koszul type} \cite{Shima}. 
This case corresponds to gluing regular models by affine Legendre equivalence of the form $\mathcal{L}_F$ with $\bc=0$ but possibly $d\not=0$. 
\end{rem}

\begin{exam}\upshape \label{exponential_family} 
 (Amari \cite{AmariNagaoka00, Amari16}). 
An {\em exponential family} $M$ is a family of probability density functions of the form 
$$p(\bu|\theta)=\exp(\bu^T\theta - \psi(\theta))$$ 
where $\bu=(u_1, \cdots, u_n) \in \R^n$ is a random valuable (with its measure $d\mu$) 
and $\theta=(\theta_1, \cdots, \theta_n)^T \in U\subset \R^n$ are parameters ($U$ is an open set). 
The normalization factor $\psi(\theta)=\log \int \exp(\bu^T\theta)\, d\mu$ is called the potential of this family. Fix the affine structure of $U$, and 
put $\rd_i=\frac{\rd}{\rd \theta_i}$. We see that the expectation is the corresponding dual coordinate 
$$\eta_i:={\bf E}[u_i| \theta]=\rd_i\psi(\theta)$$
and the (co)variance are written by 
$$ h_{ij}:={\bf V}[\bu| \theta]_{ij}=\rd_i\rd_j \psi(\theta)={\bf E}\left[(\rd_i\log p)(\rd_j\log p)\right],$$
where the last one means the {\em Fisher-Rao information}. 
If $h=[h_{ij}]$ is positive and one regards $\theta, \eta$ as the $\nabla, \nabla^*$-affine coordinates, respectively, 
then $(M, h, \nabla, \nabla^*)$ becomes a dually flat manifold. 
Normal distributions and finite discrete distributions are typical examples. 
\end{exam}

\section{Quasi-Hessian structure}\label{sec:3}

Our main idea is to consider not only regular models but also general Legendre submanifolds of $\R^{2n+1}$. 
Then the Lagrange-Legendre singularity theory naturally comes up into the picture (Arnol'd el al \cite{AGV}, Izumiya-Ishikawa \cite{Izumiya_Ishikawa}). Nevertheless, in this paper, we only use very basic notions/properties in the theory, which are prepared in \S \ref{sec:3-1}.  As another new ingredient, we introduce in \S \ref{sec:3-2} an affine geometric version of the coherent tangent bundle in Saji-Umehara-Yamada \cite{SUY}. In \S \ref{sec:3-3} and \S \ref{sec:3-4} we define a quasi-Hessian manifold endowed with a particular cubic tensor. 

\subsection{$e/m$-wavefronts and $e/m$-caustics}\label{sec:3-1}
A {\it Legendre map} is the composition 
$$\pi\circ \iota: L \to N \to B$$ 
of the inclusion $\iota$ of a Legendre submanifold $L$ and the projection of a Legendre fibration $\pi: N \to B$. 
The image is usually called a {\it wavefront}; we denote it by $W(L)$ in this paper. 
The Legendre map $\pi\circ \iota: L \to B$ may have singular points, 
i.e., points $p$ on $L$ at which the rank of the differential is not maximum (equivalently, $T_pL$ is tangent to the fiber of $\pi$), called {\em Legendre singularities} \cite{AGV, Izumiya_Ishikawa}. 
Then the wavefront is no longer a submanifold. 

From now on, we consider the diagram $(dL)$ of double Legendre fibrations on $\R^{2n+1}$ and an {\em arbitrary} Legendre submanifold $L \subset \R^{2n+1}$. 
So we have two Legendre maps
$$\pi^e:=\pi\circ \iota: L \to \R^n_{\bbx}\times \R_z, \quad \pi^m:=\pi'\circ \iota:  L \to \R^n_{\bbp}\times \R_{z'}$$ 
and call them the {\em $e$- and $m$-Legendre maps}, respectively, 
following a traditional notation in information geometry  
(``$e$-" and ``$m$-" come from words in statisitcs, i.e., exponential and mixture) \cite{AmariNagaoka00}. 

\begin{dfn}\upshape
{\bf ($e/m$-wavefronts)} 
We set 
$$W_e(L):=\pi^e(L)\subset \R^n_{\bbx}\times \R_z, \;\;\; 
W_m(L):=\pi^m(L)\subset \R^n_{\bbp}\times \R_{z'}, $$
and call them the {\em $e/m$-wavefronts associated to $L$}, respectively. 
\end{dfn}

The $e/m$-wavefronts are {\em Legendre dual} to each other in point-hyperplane duality principle (Remark \ref{proj}). 

Usually,  the projection of a Lagrange submanifold of $T^*\R^n$ to the base is called a {\em Lagrange map} \cite{AGV, Izumiya_Ishikawa}. 
So we have the {\em $e/m$-Lagrange maps} 
$$\pi^e_1=\Pi\circ \pi^e: L \to \R^n_{\bbx}, \quad \pi^m_1=\Pi\circ \pi^m: L \to \R^n_{\bbp}.$$  
It is easy to see that the following two conditions on points $p\in L$ are equivalent: 
\begin{itemize}
\item  $p$ is a singular point of the Legendre map $\pi^e: L \to \R^n_{\bbx}\times \R_z$; 
\item $p$ is a singular point of the Lagrange map $\pi^e_1: L \to \R^n_{\bbx}$. 
\end{itemize}
Indeed,  any $v \in T_pL$ enjoys $dz_p(v)-\bp(p)^T d\bx_p(v)=0$, thus, 
if $d\bx_p(v)=0$, then $dz_p(v)=0$.

\begin{dfn}\upshape
{\bf ($e/m$-caustics)}
The {\em $e$-critical set} $C(\pi^e_1)\subset L$ consists of all singular points of the $e$-Lagrange map $\pi^e_1: L \to \R^n_{\bbx}$, 
and we call its image $\pi^e_1C(\pi^e_1)\subset \R^n_{\bbx}$ the {\em $e$-caustics associated to $L$}. 
The $m$-version is defined in entirely the same way. 
\end{dfn}

\begin{dfn}\upshape
We say that $L$ is {\em locally a regular model around $p \in L$} if there is an open neighborhood of $p$ in $L$ which is a regular model of $\R^{2n+1}$, i.e., $p$ is neither {\em $e$-critical} nor {\em $m$-critical}. 
\end{dfn}

Consider the case that $p \in L$ is not $e$-critical but $m$-critical (see toy examples in Examples \ref{A2}, \ref{A3} below). 
Then $W_e(L)$ is the graph of some local potential function $z=f(\bx)$ defined near $\pi^e_1(p)$. 
Take $\bx$ as local coordinates of $L$ around $p$. Then the $e$-Lagrange map is written as the identity map of $\bx$ and the $e$-caustics is empty, while the $m$-Lagrange map $\pi^m: L \to \R^n_{\bbp}$ is written 
as the gradient map $\bp=\nabla f(\bx)$. Now it is critical at $p$, so $W_m(L)$ is singular. 
In this case we call $L$ a {\em  model with degenerate potential}. 
In particular, if $f$ admits inflection points in strict sense, $W_m(L)$ is the graph of a {\em multi-valued function} $z'=\varphi(\bp)$ branched along the $m$-caustics in $\R^n_{\bbp}$. 

\begin{figure}
 \includegraphics[width=10.5cm, pagebox=cropbox]{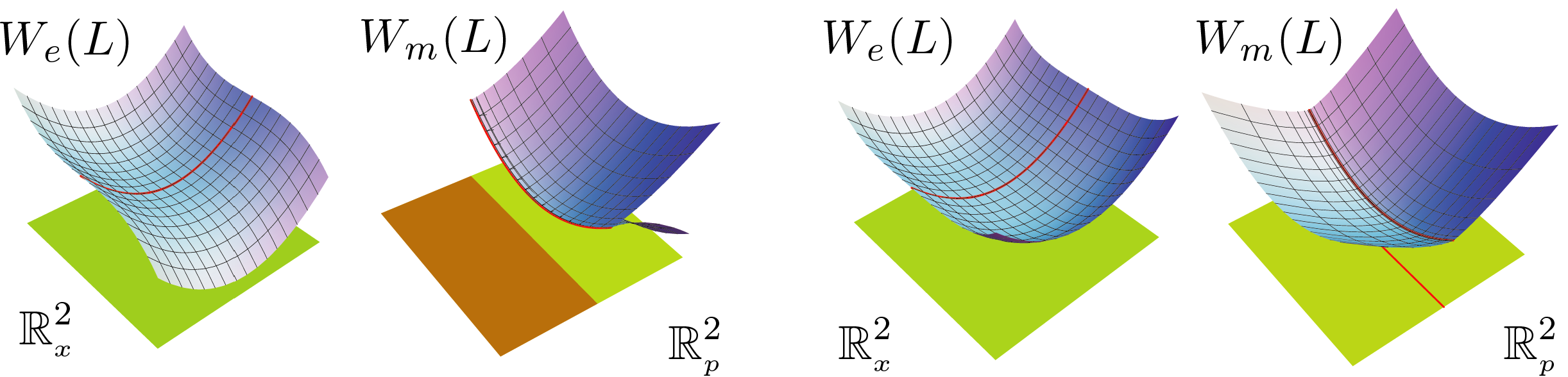}\\
 \caption{\small  
The $e/m$-wavefronts and the $e/m$-caustics (Examples \ref{A2} and \ref{A3}). 
}\label{fig1}
\end{figure}

\begin{exam}\label{A2}\upshape
{\bf ($A_2$-singularity).} Let 
$$f(\bx) = \frac{x_1^3}{3} + \frac{x_2^2}{2}.$$ 
Then $\bp=(p_1, p_2)=(x_1^2, x_2)$ and the degeneracy locus $\Sigma$ of the Hessian $h=\nabla^2f$ is defined by $x_1=0$. See the pictures on the left in Fig.\,\ref{fig1}. 
\begin{itemize}
\item[-] The {\bf $e$-wavefront} $W_e(L)$ is smooth and has parabolic points along $\Sigma$. There is no {\bf $e$-caustic}. 
 \item[-] The {\bf $m$-wavefront} $W_m(L)$ is a singular surface with {\em cuspidal edge}; it is the graph of the bi-valued dual potential 
 $$z'=\bp^T\bx-z=\frac{2}{3}x_1^3+\frac{1}{2}x_2^2=\pm\frac{2}{3}p_1^{3/2}+\frac{1}{2}p_2^2$$ 
 defined on $p_1\ge 0$ and branched along the {\bf $m$-caustics} $p_1=0$. 
\end{itemize}
This singularity does not appear, if the Hessian is be non-negative. Note that for every statistical model, the Fisher-Rao metric is non-negative.  
\end{exam}

\begin{exam}\label{A3}\upshape
{\bf  ($A_3$-singularity).} Let  
$$f(\bx) = \frac{x_1^4}{4} + \frac{x_2^2}{2}.$$ 
Then $\bp=(p_1, p_2)=(x_1^3, x_2)$ and the degeneracy locus $\Sigma$ of the Hessian  $h=\nabla^2f$ is defined by $x_1=0$. See the pictures on the right in Fig.\,\ref{fig1}. 
\begin{itemize}
\item[-] The {\bf $e$-wavefront} $W_e(L)$ is smooth and convex. There is no {\bf $e$-caustic}. 
\item[-] The {\bf $m$-wavefront}  $W_m(L)$ is a singular surface; it is the graph of the dual potential 
$$z'=\bp^T\bx-z=\frac{3}{4}x_1^4+\frac{1}{2}x_2^2=\frac{3}{4}p_1^{4/3}+\frac{1}{2}p_2^2,$$ 
which is defined on the entire space but singular along the {\bf $m$-caustics} $p_1=0$. 
\end{itemize}
This is a typically degenerate minimum of functions 
and also a typical type of singularities with $\Z_2$-symmetry (cf. \cite{AGV}). 
\end{exam}

Furthermore, it can happen that $p \in L$ is $e$-critical and $m$-critical simultaneously. 
Then both wavefronts $W_e(L)$ and $W_m(L)$ become singular at $\pi^e(p)$ and $\pi^m(p)$, respectively. 
In general, by Implicit Function Theorem, there is a partition $I \sqcup J=\{1, \cdots, n\}$ ($I\cap J=\emptyset$) such that $L$ is locally parametrized  around $p$ by coordinates $\bx_I=(x_i)^T$ and $\bp_J=(p_j)^T$ ($i \in I$, $j \in J$). 
In fact, we can find a function $g(\bx_I, \bp_J)$ such that  $L$ near $p$ is expressed by 
$$\bp_I= \frac{\rd g}{\rd \bx_I}, \;\; \bx_J=- \frac{\rd g}{\rd \bp_J}, \;\; 
z=\bp_J^T\bx_J+g(\bx_I, \bp_J),$$
where we write $\frac{\rd g}{\rd \bbx_I}=(\frac{\rd g}{\rd x_{i_1}}, \cdots)^T \; (I=(i_1, i_2, \cdots))$. 
This follows from the form (1) in \S \ref{sec:2-1} and  
the {\em canonical transformation} 
$$\R^{2n+1} \to \R^{2n+1} \;\; 
(\bx, \bp, z) \mapsto  (\bx_I, \bp_J, \bp_I, -\bx_J, -\bp_J^T\bx_J+z)$$
which preserves the contact structure. 
Usually, $g(\bx_I, \bp_J)$ is called 
a {\em generating function of $L$ around $p$}  \cite[\S 20]{AGV}. 
In particular, in case that 
$J=\emptyset$ (resp. $I=\emptyset$), 
a generating function is a {\em potential} $z=f(\bx)$ (resp. {\em dual potential}  $z'=\varphi(\bp)$). 
The $e/m$-Legendre maps are locally expressed as follows. 
$$\begin{array}{ll}
\pi^e: (L, p) \to \R^n_{\bbx}\times \R_z,  &
\displaystyle  
(\bx_I, \bp_J) \mapsto 
\left(\bx_I, \; - \frac{\rd g}{\rd \bp_J}, \; - \bp_J^T\frac{\rd g}{\rd \bp_J}+g(\bx_I, \bp_J)\right), \\
\pi^m: (L, p) \to \R^n_{\bbp}\times \R_{z'},  &
\displaystyle  
(\bx_I, \bp_J) \mapsto 
\left(\frac{\rd g}{\rd \bx_I},\; \bp_J, \; \bx_I^T\frac{\rd g}{\rd \bx_I}-g(\bx_I, \bp_J)\right).
\end{array}$$
Also the $e/m$-Lagrange maps $\pi^e_1, \pi^m_1$ are obtained by ignoring the last $z$ and $z'$-coordinate, respectively.

\begin{exam}\label{AA}\upshape
Let  
$$g(x_1, p_2) = \frac{x_1^3}{3} + \frac{p_2^4}{4}$$ 
be a generating function. 
The $e/m$-Legendre maps $\pi^e$ and $\pi^m$ send $(x_1, p_2)$ to 
$$(x_1, x_2, z)=\left(x_1, -p_2^3, \frac{x_1^3}{3}-\frac{3p_2^4}{4}\right), 
\quad (p_1, p_2, z')=\left(x_1^2, p_2, \frac{2x_1^3}{3}-\frac{p_2^4}{4}\right),$$
respectively, so those images $W_e(L)$ and $W_m(L)$ are singular surfaces having some own geometric nature, and the $e/m$-caustics are defined by $x_2=0$ on $\R^2_{\bbx}$ and $p_1=0$ on $\R^2_{\bbp}$, see Fig.\,\ref{fig11}. 
\begin{figure}[h]
 \includegraphics[width=7cm, pagebox=cropbox]{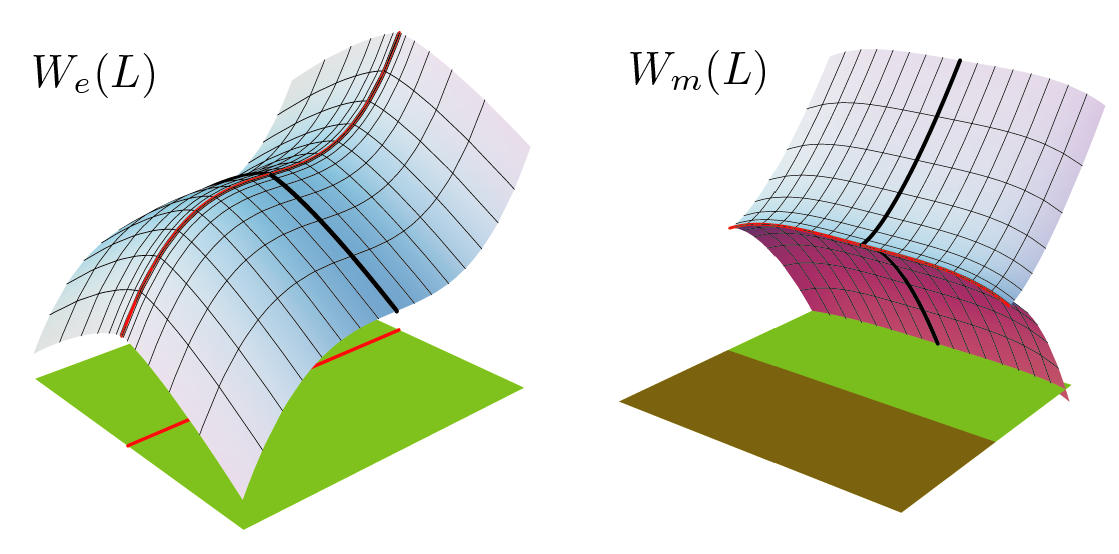}\\
 \caption{\small 
Both $e/m$-wavefronts are singular (Example \ref{AA}). 
}\label{fig11}
\end{figure}
\end{exam}

\begin{rem}\upshape
\label{hierarchy}
{\bf (Hierarchical structure)} 
For a dually flat manifold $L$ with a potential $z=f(\bx)$ (i.e., a regular model), 
there are two systems of coordinates, $\bx$ and $\bp\, (=\frac{\rd f}{\rd \bbx})$, 
which are $\nabla$-flat and $\nabla^*$-flat, respectively. 
That produces a {\em hierarchical structure} -- we may take a new system of coordinates $(\bx_I, \bp_J)$,  called  {\em mixed coordinates} in Amari \cite{Amari16}, which yields two foliations of complementary dimensions on $L$ defined by $\bx_I=const.$ and $\bp_J=const.$; their leaves are $\nabla^*$-flat and $\nabla$-flat, respectively, and mutually orthogonal with respect to the Hessian metric associated to $f$. 
This structure is useful for application, see  \cite{Amari16}. 
For an arbitrary Legendre submanifold $L$, a potential may not exist globally, but as seen above, for any $p \in L$, 
we can always find a generating function $g(\bx_I, \bp_J)$ on a neighborhood $U$ of $p$ (possible choices of the partition $I \sqcup J$ depends on $p$). That {\em locally} defines mixed coordinates $(\bx_I, \bp_J)$ and two orthogonal foliations on $U$ (see Remark \ref{h2} below). 
Usually these coordinates can not be extended to the entire space $L$, because of the presence of $e/m$-caustics (i.e., $h$ is degenerate). Nevertheless, this new structure is well organized globally, that we will formulate properly in the following subsections. 
\end{rem}

\subsection{Coherent tangent bundles}\label{sec:3-2}

Let $L$ be a Legendre submanifold of $\R^{2n+1}$. 
As seen above,  the $e$-wavefront $W_e(L)$ is not a manifold in general, but there is an alternate to its `tangent bundle'. 
Every point $p \in L$ defines a hyperplane $E_p$ in $T_{\pi^e(p)}(\R^n_{\bbx}\times \R_z)=\R^n_{\bbx}\times \R_z$, and the family of such hyperplanes form a vector bundle on $L$ of rank $n$: 
$$
E (=E_L) := \{\;(p, w) \in L\times (\R^n_{\bbx} \times \R_z)\; \mid \;dz_p(w)-\bp(p)^T d\bx_p(w)=0\; \}. 
$$ 
Since $L$ is Legendrian, we see 
$$T_pL \subset \xi_p=\ker \theta_p=(d\pi_p)^{-1}(E_p),$$ 
thus 
$d\pi^e(T_pL) \subset E_p$. 
We then associate a vector bundle map (a smooth fiber-preserving map which is linear on each fiber) 
$$\Phi : TL \to E, \quad v_p \mapsto d\pi^e_p(v_p).$$ 
Note that $\Phi$ is isomorphic if and only if $\pi^e: L \to \R^n_{\bbx}\times \R$ is an immersion. 

\begin{rem}\upshape
\label{E}
We remark that $E$ is the ``limiting" tangent bundle of the $e$-wavefront $W_e(L)$. 
Note that the kernel of $d\pi_p: T_p\R^{2n+1} \to T_{\pi(p)}(\R^n_{\bbx}\times \R)$ is spanned by $\frac{\rd}{\rd p_i}$'s. 
If $p \in L$ is a regular point of the $e$-Legendre map $\pi^e: L \to \R^n_{\bbx}\times \R$, then 
$T_pL \cap \ker d\pi_p=\{0\}$ and 
$$E_p={\rm Im}\, d\pi^e_p(T_pL) =T_{\pi^e(p)}W_e(L).$$ 
In fact, in this case, $\pi^e$ is an immersion around $p$, so $W(L)$ is a submanifold around $\pi^e(p)$. 
If a sequence of regular points $p_n \in L$ of $\pi^e$ converges to a critical point $p$, then the image of $T_{p_n}L$ converges to $E_p$ (in the Grassmannian of $n$-planes in $\R^{n+1}=\R_{\bbx}^n\times \R$) because of the continuity of the bundle $E$. 
In this case, $W(L)$ is singular at $\pi^e(p)$, thus the tangent space at that point is not defined, 
but it has the limiting tangent space $E_p$ as an alternate. 
Another characterization of $E_p$ is 
$$E_p=\ker\left[d\pi'_p: T_p\R^{2n+1} \to T_{\pi^m(p)}(\R^n_{\bbp}\times \R_{z'})\, \right]$$ 
through the inclusion $E_p \subset \R_{\bbx}^n \oplus \{0\} \oplus  \R_z \subset \R^{2n+1}=T_p\R^{2n+1}$. 
In fact, the Jacobian matrix of $\pi'$ at $p$ is 
$$\left[
\begin{array}{ccc}
O&E&0\\
\bp(p)^T&\bx(p)^T&-1
\end{array}
\right]$$
so its kernel is given by $dz_p-\bp(p)^T d\bx_p=0$ and $d\bp_p=0$. 
Note that the contact hyperplane splits as 
$$\xi_p=\ker d\pi'_p\oplus \ker d\pi_p.$$ 
\end{rem}

Let $\widetilde{\nabla}$ be the flat connection on the affine space $\R^n_{\bbx}\times \R$, and 
for any $p \in L$, 
let $\psi_{p} : \mathbb{R}^n_{\bbx} \times \mathbb{R} 
\to E_{p}$ denote 
the linear projection along the $z$-axis. 
Then a connection $\nabla^E$ of the vector bundle $E$ over $L$ is naturally defined by
$$\nabla^E_X\xi (p) := \psi_{p} \circ  \widetilde{\nabla}_X\xi(p)$$ 
where $\xi$ is any section of $E$ and $X$ is any vector field on $L$ around $p$. 

\begin{lem}\label{nabla}
The resulting connection $\nabla^E$ is flat and `relatively torsion-free', i.e.,  
for any vector fields $X, Y$ on $L$, it holds that 
$$\nabla^E_X(\Phi(Y))-\nabla^E_Y(\Phi(X))=\Phi([X, Y]).$$
\end{lem}

\proof 
Put $s_i(p):=\frac{\rd}{\rd x_i}+p_i\frac{\rd}{\rd z} \in E_p$ ($1\le i \le n$), 
then they form a frame of flat global sections of $E$: 
$$\nabla^E_Xs_i=\psi_p(\tilde{\nabla}_X s_i)= \psi_p(X(p_i)\frac{\rd}{\rd z})=0.$$ 
Thus $\nabla^E$ is flat. 
Next, a key point is that 
$\Phi$ is represented by the Jacobi matrix of the $e$-Lagrange map $\pi^e_1=(f_1, \cdots, f_n): L \to \R^n_{\bbx}$, 
i.e., 
$\Phi(\rd_j)=\sum_{i=1}^n (\rd_jf_i) s_i$ 
in local coordinates $(t_1, \cdots, t_n)$ of $L$ with $\rd_j=\frac{\rd}{\rd t_j}$. 
Let $X=\sum_{k=1}^n a_k \rd_k$ and $Y=\sum_{j=1}^n b_j \rd_j$, 
then 
$$\nabla^E_X \Phi(Y)=\sum_{i,j,k} ((\rd_k\rd_j f_i)a_kb_j+(\rd_jf_i)a_k(\rd_kb_j))s_i.$$ 
The rest is shown by a direct computation. 
\qed

\

\begin{dfn}\upshape 
\label{coherent}
We call $(E, \Phi, \nabla^E)$ the {\em coherent tangent bundle associated to the $e$-wavefront $W_e(L)$}. 
\end{dfn}

\begin{rem}\upshape 
The definition of coherent tangent bundles is originally due to Saji-Umehara-Yamada \cite[\S 6]{SUY} from the viewpoint of Riemannian geometry. They have studied several kinds of curvatures associated to wavefronts. In our case, we use the fixed affine structure of the ambient space of the wavefront, instead of metric. Also affine differential geometry of wavefronts should be rich. 
\end{rem}

In entirely the same way, 
for the $m$-wavefront $W_m(L)$, 
the coherent tangent bundle $E'$ with $\Phi':=d\pi^m : TL \to E'$ and $\nabla^{E'}$ is defined: 
$$
E' (=E'_L) := \{\;(p,w) \in L\times (\mathbb{R}^n_{\bbp} \times \mathbb{R}_{z'})\; \mid \;dz'(w)-\bx(p)^T d\bp(w)=0\; \}. 
$$
In fact, the double Legendre fibration $(dL)$ can be viewed as the pair of maps $(\pi\circ \mathcal{L}^{-1}, \pi)$ 
using different coordinates $(\bp, \bx, z')$ of $\R^{2n+1}$, and then the above construction yields $(E', \Phi', \nabla^{E'})$ 
in this dual side. 
In particular, $E'_p$ is identified with $\ker d\pi^e_p$ (see Remark \ref{E}). 

We have defined $E$ and $E'$ as vector bundles on $L$, although they are actually defined on the ambient space $\R^{2n+1}$. 
The contact hyperplane $\xi$ has the direct sum decomposition: 
$$\xi_p=\ker d\pi'_p\oplus \ker d\pi_p \simeq E_p \oplus E'_p \simeq \R^n_{\bbx}\oplus \R^n_{\bbp}.$$
Here we have canonical frames of flat sections for both $E$ and $E'$, 
$$s_i(p) = \frac{\rd}{\rd x_i} + p_i\frac{\rd}{\rd z} \in E_p,\quad s_i^*(p)= \frac{\rd}{\rd p_i}+x_i\frac{\rd}{\rd z'} \in E'_p,$$ 
by which $E'$ is identified with the dual to $E$ and vice-vasa, and 
there are natural correspondence 
$s_l \leftrightarrow \frac{\rd}{\rd x_l}$ and $s_l^* \leftrightarrow \frac{\rd}{\rd p_l}$ via projections along the $z$ and $z'$-axes. 
The vector bundle $\xi$ carries not only the symplectic form 
$\omega=\sum_{i=1}^n dx_i\wedge dp_i$ 
but also a pseudo-Riemannian metric of type $(n,n)$ 
induced from 
$$\tau:=\sum_{i=1}^n dx_idp_i=\frac{1}{2}\sum_{i=1}^n (dx_i \otimes dp_i+dp_i \otimes dx_i).$$
Using frames $s_i$ and $s_j^*$, 
we may write vectors of $E_p$ and $E'_p$ as column vectors $\bu$ and $\bu'$, respectively, and then 
$$\tau(\bu\oplus \bu', \bv\oplus \bv')=
\frac{1}{2}
\left[\, \bu^T \; \bu'^T \right] 
\left[
\begin{array}{cc}
O&E\\
E&O
\end{array}
\right] 
\left[
\begin{array}{cc}
\bv\\
\bv'
\end{array}
\right] 
=\frac{1}{2}(\bu^T\bv'+\bu'^T\bv)
$$
and also $\omega(\bu\oplus \bu', \bv\oplus \bv')=\frac{1}{2}(\bu^T\bv'- \bu'^T\bv)$. 

Any affine Legendre equivalence $\mathcal{L}_F$ preserves $\omega$ and $\tau$ on  $\xi$, 
because it sends $\bu\oplus \bu'$ to $A\bu\oplus A'\bu'$ with $A'=(A^T)^{-1}$. 

\begin{dfn}\upshape
We define the {\em quasi-Hessian metric} of $L$ by the pullback of $\tau$: 
$$h(Y, Z):=\tau(\iota_*Y, \iota_*Z)\;\;\;  \mbox{for $Y, Z \in TL$} $$
where $\iota_*=\Phi\oplus \Phi': TL \hookrightarrow \xi=E\oplus E'$ is the inclusion (it is a Lagrange subbundle). 
\end{dfn}

Note that $h$ is a possibly degenerate symmetric $(0,2)$-tensor, although we abuse the word `metric'. 
If $\Phi$ is isomorphic, then $h$ exactly coincides with the Hessian metric associated to a potential $z=f(\bx)$; 
any vector of $T_pL$ is written by $\bu\oplus H\bu \in \xi_p$ where $H=[h_{ij}]$ is the Hessian matrix, thus 
$$h(\bu, \bv)=\tau(\bu\oplus H\bu, \bv\oplus H\bv)=\bu^TH\bv.$$ 
In general,  a local expression of $h$ is given as follows. 

\begin{lem}\label{h}
Let $g(\bx_I, \bp_J)$ be a generating function. Then, 
$$h=\sum_{i, k\in I}\frac{\rd^2 g}{\rd x_i\rd x_{k}}\,dx_idx_{k} - \sum_{j, l\in J}\frac{\rd^2 g}{\rd p_j\rd p_{l}}\,dp_jdp_{l}.$$
\end{lem}

\proof 
A direct computation shows 
\begin{eqnarray*}
\tau &=& d\bx_I^Td\bp_I+d\bx_J^Td\bp_J\\
&=&d\bx_I^Td(\rd_Ig)-d(\rd_J g)^Td\bp_J\\
&=&d\bx_I^Tg_{II}d\bx_I+d\bx_I^Tg_{IJ}d\bp_J
-(g_{JI}d\bx_I)^Td\bp_J-(g_{JJ}d\bp_J)^Td\bp_J\\
&=&d\bx_I^Tg_{II}d\bx_I-d\bp_J^Tg_{JJ}d\bp_J.
\end{eqnarray*}
Here we use the notation of symmetric products of $1$-forms and $(g_{JI})^T=g_{IJ}$. 
\qed

\

\begin{lem}\label{g}
Let $p \in L$. The following properties are equivalent: 
\begin{enumerate}
\item[(1)] $h$ is non-degenerate at $p$ ; 
\item[(2)] $p$ is neither of $e$-critical nor $m$-critical; 
\item[(3)] $L$ is locally a regular model around $p$; 
\item[(4)] $h$ is the Hessian metric associated to a local potential $z=f(\bx)$ near $p$; 
\item[(5)] both $\Phi$ and $\Phi'$ are isomorphisms at $p$. 
\end{enumerate}
\end{lem}

\proof By Lemma \ref{h},  (1) means that both $g_{II}$ and $g_{JJ}$ are non-degenerate. Then, 
using normal forms of the $e/m$-Lagrange maps $\pi^e_1$ and $\pi^m_1$ written in the end of \S 2.3, 
those maps are locally diffeomorphic by Inverse Mapping Theorem, so it is just (2) and (3).  
That means that we can take a local potential $z=f(\bx)$ as generating function, that is equivalent to (4). 
Since $\Phi$ and $\Phi'$ are expressed by the Jacobi matrices of the $e/m$-Lagrange maps, 
(2) and (5) are the same. 
\qed

\begin{rem}\upshape
\label{h2}
As noted in Remark \ref{hierarchy}, locally we always find mixed coordinates $(\bx_I, \bp_J)$. 
By Lemma \ref{h}, even if $h$ is degenerate, 
leaves $\bp_J=const.$ and $\bx_I=const.$ are orthogonal: $h(\frac{\rd}{\rd x_i}, \frac{\rd}{\rd p_j})=0$ 
($i \in I$, $j\in J$).  
\end{rem}

\begin{dfn}\upshape
Let $\Sigma\, (=\Sigma_{L, h})$ denote the set of $p \in L$ at which $h$ is degenerate, equivalently, 
the locus where either $\Phi$ or $\Phi'$ is not isomorphic: 
$$\Sigma=C(\pi^e) \cup C(\pi^m).$$ 
We call $\Sigma$ the {\em degeneracy locus} of the quasi-Hessian metric $h$. 
\end{dfn}

Since $L$ is Legendrian, $T_pL$ is a Lagrange subspace of the symplectic vector space $\xi_p= E_p \oplus E'_p$.  
Note that $\Phi$ (resp. $\Phi'$) is the linear projection of $T_pL$ to the factor $E_p$ (resp. $E_p'$), 
and especially, $\dim T_pL\cap E_p\ge 1$  (resp. $\dim T_pL\cap E'_p\ge 1$) if and only if 
$p$ is $m$-critical  (resp. $e$-critical).  In particular,  the null space of $h$ splits: 
$${\rm null}\, h_p=\ker \Phi'_p \oplus \ker \Phi_p=(T_pL\cap E_p) \oplus (T_pL\cap E'_p).$$ 

\begin{dfn}\upshape \label{duallyflat2} 
For an {\em arbitrary} Legendre submanifold $L \subset \R^{2n+1}$,  
we call the triplet $(h, (E, \nabla^E, \Phi),  (E', \nabla^{E'}, \Phi'))$ the {\em dually flat structure} of $L$. 
$$
\xymatrix{
&TL \ar[dl]_{\Phi} \ar[dr]^{\Phi'} &\\
E_L&&E_L'
}
$$
\end{dfn}

\begin{rem}\upshape 
\label{dually_flat_rem}
Given a regular model $L$, we have the triple $(h, E, E')$, where 
both $\Phi$ and $\Phi'$ are isomorphic. 
That restores the dually flat structure in the original form (Definition \ref{duallyflat}); 
Indeed, $\nabla$ and $\nabla^*$ on $TL$ are uniquely determined by 
$$\Phi(\nabla_X Y)=\nabla^E_X \Phi(Y),  \quad 
\Phi'(\nabla^*_X Y)=\nabla^{E'}_X \Phi'(Y) $$
where $X, Y$ are arbitrary vector fields on $L$. 
On the other hand, a singular model $L$ with degenerate potential $z=f(\bx)$ is 
the case that $\Phi$ is isomorphic and $\Phi'$ is not. 
Then the connection $\nabla$ of $TL$ is obtained from $\nabla^E$ via $\Phi$ in the same way as above, 
while $\nabla^*$ does not exist. 
If both $\Phi$ and $\Phi'$ are not isomorphic, there is no connection on $TL$. 
\end{rem}

\subsection{Quasi-Hessian manifolds}\label{sec:3-3}
Our generalized dually flat structure presented in Definition \ref{duallyflat2} is compatible with affine Legendre equivalence. 
That means that if an affine Legendre equivalence $\mathcal{L}_F$ identifies Legendre submanifolds $L_1$ and $L_2$, then the quasi-Hessian metrics are preserved, $\mathcal{L}_F^*h_2=h_1$, and $\mathcal{L}_F$ naturally induces vector bundle isomorphisms between coherent tangent bundles, $E_{L_1} \simeq E_{L_2}$ and $E'_{L_1} \simeq E'_{L_2}$, such that the isomorphisms identify equipped affine flat connections and we have the following commutative diagram 
$$
\xymatrix{
E_{L_1} \ar[d]_{\simeq}&TL_1 \ar[d]^{\simeq}_{\mathcal{L_F}} \ar[l]_{\Phi_1}\ar[r]^{\Phi'_1} & E'_{L_1} \ar[d]^{\simeq} \\
E_{L_2} &TL_2 \ar[l]_{\Phi_2}\ar[r]^{\Phi'_2} & E'_{L_2} 
}
$$

Thus the ordinary gluing construction works. 
To be precise, suppose that we are given a collection $\{L_\alpha\}_{\alpha\in \Lambda}$, 
where $\Lambda$ is a countable set, such that it satisfies the following properties: 
\begin{enumerate}
\item[(i)] for every $\alpha \in \Lambda$, $L_\alpha$ itself is an open manifold and it is embedded in $\R^{2n+1}$ as a Legendre submanifold, called a {\em local model};
\item[(ii)] for every $\alpha, \beta \in \Lambda$, there is an open subset $L_{\alpha\beta} \subset L_\alpha$ (also $L_{\beta\alpha} \subset L_\beta$) and a diffeomorphism $\mathcal{L}_\alpha^\beta: L_{\alpha\beta} \to L_{\beta\alpha}$ such that 
over each connected component of $L_{\alpha\beta}$, $\mathcal{L}_\alpha^\beta$ is given by an affine Legendre equivalence of the ambient space $\R^{2n+1}$; 
\item[(iii)] for $\alpha, \beta, \gamma \in \Lambda$, it holds that $\mathcal{L}_\alpha^\gamma=\mathcal{L}_\beta^\gamma\circ \mathcal{L}_\alpha^\beta$ on $L_{\alpha\beta} \cap L_{\alpha\gamma}$. 
\end{enumerate}
Let $M$ be the resulting topological space from these gluing data $\mathcal{U}=\{L_\alpha, \mathcal{L}_\alpha^\beta\}$. Assume that $M$ is Hausdorff, then $M$ itself becomes an $n$-dimensional manifold in the ordinary sense. 
One can naturally associate a possibly degenerate $(0,2)$-tensor $h$ on $M$ and a pair of globally defined dual coherent tangent  bundles $E$ and $E'$ on $M$ with bundle maps $\Phi: TM \to E$ and $\Phi: TM \to E'$ equipped with affine flat connections.  The bundles $E$ and $E'$ are dual to each other. 

\begin{dfn}\upshape\label{quasi-Hessian_manifold}
We call $(M, \mathcal{U})$ equipped with $(h, (E, \nabla^{E}, \Phi), (E', \nabla^{E'}, \Phi'))$  a {\em quasi-Hessian manifold}. 
We define the {\em degeneracy locus} $\Sigma$ to be the locus of points of $M$ at which $h$ is degenerate. 
\end{dfn}

Since the gluing maps $\mathcal{L}_\alpha^\beta$ also act on a neighborhood of $L_{\alpha\beta}$ in $\R^{2n+1}$  and preserve the contact structure, 
$M$ is realized as a Legendre submanifold in some ambient contact manifold. 

By the above construction, it is obvious to see  

\begin{prop}\upshape 
Let $M$ be a quasi-Hessian manifold. Then $h$ is non-degenerate everywhere 
if and only if $M$ is a Hessian manifold. 
\end{prop}

\proof 
From the equivalence of (1) and (3) in Lemma \ref{g}, we see that  
$h$ is non-degenerate everywhere if and only if any local models $L_\alpha$ are regular models, 
that means $M$ is a Hessian manifold (Proposition \ref{Hessian_manifold}). \qed

\begin{rem}\label{singular_L}
\upshape 
More generally, we may allow a local model $L_\alpha$ not to be a manifold  but a {\em singular Legendre variety}; it is a closed subset with a partition (stratification) into integral submanifolds of the contact structure  (the projection to the cotangent bundle is called a {\em singular Lagrange variety}),  see, e.g., Ishikawa \cite{Ishikawa92}. That results a quasi-Hessian manifold with singularities. 
\end{rem}

An intrinsic definition of quasi-Hessian manifolds is also available. Roughly speaking, it is an $n$-dimensional manifold $M$ equipped with a pair of flat coherent tangent bundles $(E, \nabla^E, \Phi)$ and $(E', \nabla^{E'}, \Phi')$ of rank $n$; we impose two conditions: 
\begin{enumerate}
\item[(a)]
 the vector bundle $E\oplus E'$ of rank $2n$ is endowed with a symplectic structure $\omega$ and a pseudo-Riemannian metric $\tau$ of type $(n,n)$ satisfying $\omega(u,v)=\tau(u,v)=0$ ($u, v\in E_p$ or $E'_p$) and 
$\omega(u,v)=\tau(u,v)$ ($u\in E_p$, $v \in E'_p$); this condtion defines the dualily between $E$ and $E'$;  
\item[(b)] the bundle map 
$$\Phi\oplus \Phi': TM \to E \oplus E'$$ 
is injective and the image is a Lagrange subbundle which is certainly {\em integrable} in order to ensure to find a local model around each point of $M$ as in Definition \ref{quasi-Hessian_manifold}. We omit the detail here. 
\end{enumerate}
This also suggests a degenerate version of the so-called Codazzi structure (cf. \cite{Shima}).

\subsection{Cubic tensor and $\alpha$-family}\label{sec:3-4}
In the theory of dually flat manifolds \cite{AmariNagaoka00}, 
not only the Hessian metric $h$ but also the Amari-Chentsov tensor $T:=\nabla h$ takes an essential role; it satisfies 
$$T(X, Y, Z)=h(\nabla^*_XY, Z)-h(\nabla_XY, Z)=h(Y, \nabla^*_XZ)-h(Y, \nabla_XZ)$$ 
for vector fields $X, Y, Z$. 

Note that whenever $\nabla$ exists, 
the tensor $T$ is defined everywhere, independently whether or not $h$ is non-degenerate. 
This is an easy case. We generalize the Amari-Chentsov tensor for an arbitrary quasi-Hessian manifold 
$$(M, h, (E, \nabla^E, \Phi), (E', \nabla^{E'}, \Phi'))$$
but the way is not obvious at all, because there is no connection of $TM$. 
Finally we will see that the obtained tensor is a very natural one (Proposition \ref{C_local} below). 

\begin{lem}\label{tau}
For any vector field $X$ on $M$, and for any sections $\eta$ of $E$ and $\zeta'$ of $E'$, 
it holds that 
$$X\tau(\eta, \zeta')=\tau(\nabla^E_X\eta, \zeta')+ \tau(\eta, \nabla^{E'}_X\zeta')$$
where we put $\tau(\eta, \zeta'):=\tau(\eta\oplus 0, 0\oplus \zeta')$ for short. 
\end{lem}

\proof Take local frames of flat sections $s_i$ of $E$ and $s_j^*$ of $E'$ with $\tau(s_i, s_j^*)=\frac{1}{2}\delta_{ij}$ ($1\le i, j \le n$) on an open set $U \subset M$. 
Put $\eta=\sum a_is_i$ and $\zeta'=\sum b_js_j^*$ where $a_i, b_j$ are functions on $U$, 
then 
$$\nabla^E_X\eta=\sum X(a_i)s_i, \;\; \nabla^{E'}_X\zeta'=\sum X(b_j)s_j^*, \;\; \tau(\eta, \zeta')=\frac{1}{2}\sum a_ib_i.$$ 
This leads to the equality. 
\qed

\

For $Y, Z \in TM$, put 
$$\eta=\Phi(Y), \;\; \zeta=\Phi(Z) \; \in E, \;\; \eta'=\Phi'(Y), \;\;  \zeta'=\Phi'(Z)\; \in E'.$$ 
Then 
$$
h(Y, Z)=\tau(\eta\oplus \eta', \zeta\oplus \zeta')=\tau(\eta, \zeta')+\tau(\zeta, \eta'). 
$$
Using Lemma \ref{tau}, for vector fields $X, Y, Z$ on $M$, 
\begin{eqnarray*}
Xh(Y, Z)&=&X(\tau(\eta, \zeta'))+X(\tau(\zeta, \eta'))\\
&=&\tau(\nabla^E_X\eta, \zeta') + \tau(\eta, \nabla^{E'}_X\zeta')+ \tau(\nabla^E_X\zeta, \eta') + \tau(\zeta, \nabla^{E'}_X\eta').
\end{eqnarray*}
We call the sum of first and third terms the $\nabla^E$-part,  the rest the $\nabla^{E'}$-part, tentatively. 
We are concerned with their difference. 

\begin{dfn}\upshape 
\label{C}
For a quasi-Hessian manifold $M$, 
we define the {\em canonical cubic tensor} $C$ by the following $(0,3)$-tensor on $M$: 
\begin{eqnarray*}
&&C(X, Y, Z):=\tau(\eta, \nabla^{E'}_X\zeta')+\tau(\zeta, \nabla^{E'}_X\eta')
-\tau(\nabla^E_X\eta, \zeta')-\tau(\nabla^E_X\zeta, \eta'). 
\end{eqnarray*}
\end{dfn}

In particular, if $h$ is non-degenerate, then $\Phi(\nabla_X Y)=\nabla^E_X \Phi(Y)$ (Remark \ref{dually_flat_rem}) 
and we have 
$$
\tau(\nabla^E_X\eta, \zeta')=\tau(\Phi(\nabla_X Y), \Phi'(Z))=\frac{1}{2}h(\nabla_X Y, Z)
$$
and so on, thus it follows that the $\nabla^E$-part and the $\nabla^{E'}$-part 
are equal to, respectively, 
$$\frac{1}{2}(h(\nabla_XY, Z)+h(Y, \nabla_XZ)), \quad \frac{1}{2}(h(\nabla^*_XY, Z)+h(Y, \nabla^*_XZ)).$$ 
Hence, we see that $C$ coincides with the Amari-Chentsov tensor $T$: 
\begin{eqnarray*}
C(X, Y, Z)&=&\textstyle 
\frac{1}{2}(h(\nabla^*_XY, Z)-h(\nabla_XY, Z))+\frac{1}{2}(h(Y, \nabla^*_XZ)-h(Y, \nabla_XZ))\\
&=&T(X, Y, Z). 
\end{eqnarray*}

Using local coordinates, we write down the tensor $C$ explicitly as follows. 
Take a local model $L \subset \R^{2n+1}$ and $p \in L$. As mentioned before, 
locally around $p$, $L$ is parameterized by some local coordinates $\bx_I, \bp_J$ with a generating function $g(\bx_I, \bp_J)$.  
For the simplicity, for each $1\le k \le n$, we set 
$$\rd_k:=\frac{\rd}{\rd x_k}\; (k\in I) \;\;\;\mbox{or}\;\;\; \frac{\rd}{\rd p_k}\; (k \in J).$$

\begin{prop} \label{C_local} 
The canonical cubic tensor $C$ is locally the third partial derivative of a generating function: for any $k, l, m$, 
$$C(\rd_k, \rd_l, \rd_m)=\rd_k\rd_l\rd_m g.$$ 
In particular, $C$ is symmetric. 
\end{prop}

\proof 
This is shown by direct computation. 
The generating function yields 
a Lagrange embedding $L \to T^*\R^n$ given by 
$$\iota: (\bx_I, \bp_J)\mapsto (\bx_I, \bx_J, \bp_I, \bp_J):=(\bx_I, -\rd_J g, \rd_I g, \bp_J),$$ 
thus the differential $\iota_*: T_pL \to T_p(T^*\R^n)=\R^{n}_{\bbx}\oplus \R^n_{\bbp}$ is written as  
$$\iota_*(\rd_k)=\rd_k-\sum_{j\in J} (\rd_k\rd_j g)\frac{\rd}{\rd x_j} + \sum_{i\in I} (\rd_k\rd_i g) \frac{\rd}{\rd p_i}.$$  
Let $s_i, s_i^*\; (1\le i \le n)$ be flat sections of $E$ and $E'$ as before; $\tau(s_i, s_j^*)=\frac{1}{2}\delta_{ij}$. 
Then for $k \in I$, 
$$ \Phi(\iota_*\rd_k)=s_k-\sum_{j\in J} (\rd_k\rd_j g)s_j, \quad 
\Phi'(\iota_*\rd_k)=\sum_{i \in I} (\rd_k\rd_i g) s_i^*,$$
and for $k\in J$, 
$$ \Phi(\iota_*\rd_k)=- \sum_{j\in J} (\rd_k\rd_j g) s_j, \quad 
\Phi'(\iota_*\rd_k)=s_k^*+\sum_{i\in I} (\rd_k\rd_i g)s_i^*.$$

Put $\eta=\Phi(\iota_*\rd_l)$, $\eta'=\Phi'(\iota_*\rd_l)$, $\zeta=\Phi(\iota_*\rd_m)$, $\zeta'=\Phi'(\iota_*\rd_m)$, and $X=\rd_k$. 

For $l \in I$, $m \in J$ and any $k$,  we have 
\begin{eqnarray*}
&&\textstyle \tau(\eta, \nabla^{E'}_X\zeta')
=\tau(s_l-\sum_J (\rd_l\rd_jg) s_j, \sum_I (\rd_k\rd_m\rd_i g) s_i^*)=\frac{1}{2}\rd_k\rd_l\rd_m g,\\
&&\textstyle \tau(\zeta, \nabla^{E'}_X\eta') 
=\tau(-\sum_J (\rd_m\rd_j g)s_j, \sum_I (\rd_k\rd_l\rd_i g) s_i^*)=0,\\
&&\textstyle \tau(\nabla^E_X\eta, \zeta')
=\tau(-\sum_J (\rd_k\rd_l\rd_j g)s_j, s_m^*+\sum_I (\rd_m\rd_i g) s_i^*)
=-\frac{1}{2}\rd_k\rd_l\rd_m g, \\
&&\textstyle \tau(\nabla^E_X\zeta, \eta')
=\tau(-\sum_J (\rd_k\rd_m\rd_j g)s_j, \sum_I (\rd_l\rd_i g) s_i^*)=0.
\end{eqnarray*}
Thus, the $\nabla^{E'}$-part minus the $\nabla^E$-part gives $C(\rd_k, \rd_l, \rd_m)=\rd_k\rd_l\rd_m g$. 

For $l, m \in I$ and any $k$,  we have 
\begin{eqnarray*}
&&\textstyle \tau(\eta, \nabla^{E'}_X\zeta')
=\tau(s_l-\sum_J (\rd_l\rd_jg) s_j, \sum_I (\rd_k\rd_m\rd_i g) s_i^*)=\frac{1}{2}\rd_k\rd_l\rd_m g,\\
&&\textstyle \tau(\zeta, \nabla^{E'}_X\eta') 
=\tau(s_m-\sum_J (\rd_m\rd_j g)s_j, \sum_I (\rd_k\rd_l\rd_i g) s_i^*)=\frac{1}{2}\rd_k\rd_l\rd_m g,\\
&&\textstyle \tau(\nabla^E_X\eta, \zeta')
=\tau(-\sum_J (\rd_k\rd_l\rd_j g)s_j, \sum_I (\rd_m\rd_i g) s_i^*)
=0, \\
&&\textstyle \tau(\nabla^E_X\zeta, \eta')
=\tau(-\sum_J (\rd_k\rd_m\rd_j g)s_j, \sum_I (\rd_l\rd_i g) s_i^*)=0.
\end{eqnarray*}
Thus, $C(\rd_k, \rd_l, \rd_m)=\rd_k\rd_l\rd_m g$. 
The same is true for the case of $l, m \in J$. 
\qed

\begin{rem}\upshape
For a dually flat manifold with potential function $f$, 
the above proposition corresponds to a well known property
$$T(\rd_i, \rd_j, \rd_k)=\rd_i\rd_j\rd_k f,$$ 
with respect to $\nabla$-affine coordinates. 
In fact, a quasi-Hessian manifold is well characterized by using $h$ and $C$, that  will be discussed within the theory of (weak) contrast functions (see \S \ref{sec:3-4}). 
\end{rem}

As well known, 
for a dually flat manifold $M$, 
the family of {\em $\alpha$-connections} is defined by 
$$\nabla^{(\alpha)}=\frac{1+\alpha}{2}\nabla+\frac{1-\alpha}{2}\nabla^*$$
$(\alpha \in \R$). 
Namely, it deforms the Levi-Civita connection using $T$ linearly. 
When $\alpha=\pm 1$, $\nabla, \nabla^*$ are recovered. 
Both  $\nabla^{(\alpha)}$ and $\nabla^{(-\alpha)}$ are mutually dual and they form the so-called {\em $\alpha$-geometry} \cite{AmariNagaoka00, Matsuzoe}. 
For a quasi-Hessian manifold $M$, we have connections of $E$ and $E'$, but none of $TM$, 
thus there is no direct analogy to $\alpha$-geometry. 
Nevertheless, as an attempt,  we define a new $(0,3)$-tensor 
\begin{eqnarray*}
N^{(\alpha)}(X, Y, Z)&:=&
\frac{1+\alpha}{2} \left[\mbox{$\nabla^E$-part}\right] + \frac{1-\alpha}{2} \left[\mbox{$\nabla^{E'}$-part}\right] \\
&=&  \frac{1}{2}Xh(Y,Z)- \frac{\alpha}{2} C(X,Y,Z). 
\end{eqnarray*}
Obviously, $N^{(-1)}(X, Y, Z)$ is the $\nabla^{E'}$-part, $N^{(1)}(X, Y, Z)$ is the $\nabla^E$-part multiplied by $(-1)$, 
and a sort of duality holds: 
$$Xh(Y,Z)=N^{(\alpha)}(X, Y, Z)+N^{(-\alpha)}(X, Y, Z).$$
In general, $N^{(\alpha)}$ is not totally symmetric, for $Xh(Y, Z)$ is not so. 
If either $\Phi_p$ or $\Phi'_p$ is isomorphic, then we may take a possibly degenerate local (dual) potential around $p$ (i.e., $I$ or $J=\emptyset$) 
as generating function $g$; thus $h$ is written by the Hessian of the potential, 
and hence $Xh(Y, Z)$ is symmetric, and $N^{(\alpha)}$ is also. 
Furthermore, if $h$ is non-degenerate, i.e., $M$ is a dually flat manifold, 
we completely restore $\alpha$-geoemtry.

\section{Divergence}\label{sec:4}
Let $(M, h, (E, \nabla^E), (E', \nabla^{E'}))$ be a quasi-Hessian manifold throughout this section. 

\subsection{Geodesic-like curves}\label{sec:4-1}
Let  $c: I \to M$ be a curve, where $I \; (\not=\emptyset) \subset \R$ is an open interval, 
and set $\dot{c}(t):=\frac{d}{dt}c(t) \in T_{c(t)}M$, the velocity vector ($t\in I$). 

\begin{dfn}
\label{em-curve}
\upshape 
A curve $c: I \to M$ is called an {\em $m$-curve} if it is an immersion ($\dot{c}(t)\not=0$) and satisfies that at every $t \in I$,  vectors of $E'_{c(t)}$ 
$$\Phi' \circ \dot{c}(t), \;\; \nabla^{E'}_{\dot{c}}(\Phi' \circ \dot{c})(t), \;\; (\nabla^{E'}_{\dot{c}})^2(\Phi' \circ \dot{c})(t), \; \cdots $$
are not simultaneously zero and any two are linearly dependent. Also an {\em $e$-curve} is defined in the same way by replacing $\Phi'$ and $E'$ by $\Phi$ and $E$, respectively. 
\end{dfn}

\begin{figure}[h]
 \includegraphics[width=8cm, pagebox=cropbox]{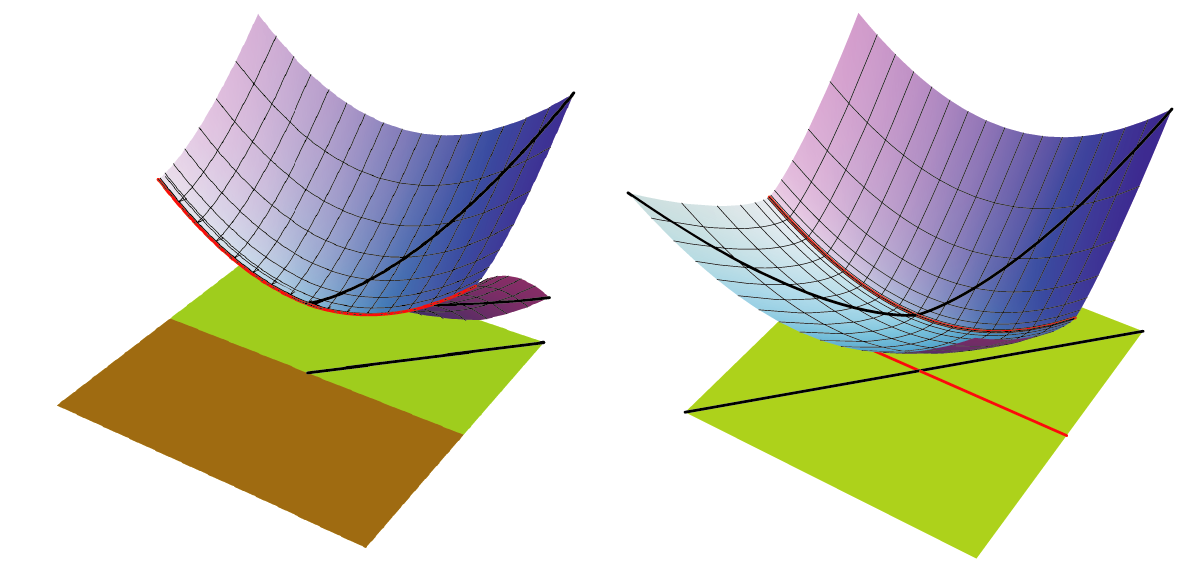}\;\; 
 \caption{\small  
The image of $m$-curves via the $m$-Lagrange map.
}\label{fig3}
\end{figure}

Suppose that the curve is given in a local model, $c: I \to L_\alpha$.  
We denote by  
$$\bp(t):=\pi^m_1\circ c(t) \in \R^n_{\bbp}$$ 
the image via the $m$-Lagrange map $\pi^m_1: L_\alpha \to \R^n_{\bbp}$. 
Note that $E'_p$ is canonically isomorphic to $\R^n_{\bbp}$ by linear projection along the $z'$-axis. 
Unless $\Phi'\circ \dot{c}(t)$ becomes to be zero,  the velocity vector $\dot{\bp}(t)$ does not vanish and its acceleration vector $\ddot{\bp}(t)$ is parallel to the velocity (it can be $0$) by the condition in Definition \ref{em-curve}. Hence $\bp(t)$ moves on a straight line in $\R^n_{\bbp}$, i.e., $c(t)$ is a {\em re-parametrization} of an $m$-geodesic (geodesic with respect to $\nabla^* = \nabla^{E'}$). 
A trouble occurs when $\Phi'\circ \dot{c}(t_0)=0$ at some $t_0$. 
Then, $\dot{\bp}(t_0)=0$, but by the condition for $m$-curve, some higher derivative is non-zero, say $\frac{d^k}{d t^k}\bp(t_0)\not=0$, and then the vector $\frac{d^{k+1}}{d t^{k+1}}\bp(t)$ is parallel to $\frac{d^k}{d t^k}\bp(t)$, so we see again that $\bp(t)$ moves on a straight line, but it meets the $m$-caustics at $t=t_0$; it stops once and then turns back or goes forward along the same line, according to $k$ even or odd,  see Fig.\,\ref{fig3} (cf. Examples \ref{A2} and \ref{A3}). 
We choose a direction vector $\bm_c$ of the straight line. For an $e$-curve $c(t)$, everything is the same, and we denote by $\be_c$ a direction of the corresponding line on $\R^n_{\bbx}$. 

\begin{rem}\upshape 
(1) Not arbitrary two points on $M$ are connected by an $m$-curve but by a {\em piecewise} $m$-curve. 
In fact, in Example \ref{A2}, you can easily find such two points on the $m$-wavefront, the left in Fig.\,\ref{fig3}. 
That is also for $e$-curves. 
(2) Take coordinates $(\bx_I, \bp_J)$ for a local model $L_\alpha$. It is easy to see that 
any $e/m$-curves satisfy a certain partial differential equation (like the geodesic equation) 
using $h=(h_{ij})$ and $C=(C_{ijk})$. 
In \S \ref{sec:3-4}, we have introduced the $\alpha$-family of cubic tensors $N^{(\alpha)}$. 
Thus we may consider an $\alpha$-analogy to $e/m$-curves; indeed, over $M-\Sigma$, it is the same as geodesics with respect to 
$\nabla^{(\alpha)}$ and $\nabla^{(-\alpha)}$. 
\end{rem}

The following definition does not depend on the choices of $L_\alpha$ and direction vectors. 

\begin{dfn} \upshape 
Let  $c_e$, $c_m$ and $\be$, $\bm$ as above. 
Let $S$ be a submanifold of $M$ and $c_m$ meets $S$ at $q \in S$. 
We say that $c_m$ and $S$ are {\em orthogonal} at $q$ if it holds that 
\begin{align*}
\bm^T d\bx(u) = 0 \;\;\; \mbox{\rm for any $u \in T_qS$}
\end{align*}
where $(\bx, \bp, z)$ is the coordinates of $\R^{2n+1}$ for a local model $L_\alpha$ containing $q$. 
Similarly, $c_e$ and $S$ are orthogonal at $q$ if $\be^T d\bp(u) = 0$ for any $u$. 
Furthermore, we say that $c_e$ and $c_m$ are {\em strictly orthogonal} at $q$ if $\be^T \bm = 0$. 
\end{dfn}

If $q \not\in \Sigma$, the above definition of the orthogonality of $c_m$ and $S$ is the same as the orthogonality with respect to the metric $h$. In fact,  taking a regular model around $q$,  the Hessian $H(q)$ is non-degenerate, and thus  
$$h(\dot{c}(t), u)=\dot{\bx}(t)^THd\bx(u)=(H\dot{\bx}(t))^Td\bx(u)=\dot{\bp}(t)^Td\bx(u)=k\bm^Td\bx(u)$$ 
for some $k\not=0$. 
However, if $q\in \Sigma$, the meaning is different in general, for it can happen that $\dot{\bp}(t_0)=0$ but $\bm\not=0$ (in this case, $\bm$ is determined by some higher derivative of $\bp(t)$). 
The reason why we define the strictly orthogonality is the same; $\be$ and $\bm$ may not be determined by velocity vectors.

\subsection{Canonical divergence}\label{sec:4-2}
Let  $L$ be a Legendre submanifold of $\R^{2n+1}$. 
We denote coordinates by 
$$p=(\bx(p), \bp(p), z(p)) \in \R^{2n+1}, \quad z'(p)=\bp(p)^T\bx(p)-z(p) \in \R.$$ 
\begin{dfn}\upshape
The {\it canonical divergence} $\mathcal{D}=\mathcal{D}_L : L \times L \to \mathbb{R}$ is defined by  
\begin{align*}
\mathcal{D}(p, q) = z(p)+z'(q) - \bx(p)^T \bp(q). 
\end{align*}
\end{dfn}
Note that $\mathcal{D}(p, p)=0$ and it is asymmetric, $\mathcal{D}(p, q)\not=\mathcal{D}(q, p)$, in general. 
In particular, if $L$ is a regular model with positve definite Hessian metric, 
this is nothing but the {\em Bregman divergence} for some convex potential $z=f(\bx)$, 
$$\mathcal{D}(p, q)=f(\bx(p))+\varphi(\bp(q)) - \bx(p)^T \bp(q),$$
where $z'=\varphi(\bp)$ is the Legendre transform of the potential \cite{AmariNagaoka00}.

\begin{lem}\label{invariance}
The canonical divergence $\mathcal{D}_L$ is invariant under affine Legendre equivalence, i.e., 
if Legendre submanifolds $L_1$ and $L_2$ of $\R^{2n+1}$ 
are affine Legendre equivalenct via $\mathcal{L}_F$, 
then it holds that 
$$\mathcal{D}_{L_1}=\mathcal{D}_{L_2}\circ (\mathcal{L}_F\times \mathcal{L}_F) \quad 
\mbox{on $L_1\times L_1$.} $$ 
\end{lem}

\proof 
Suppose that $\mathcal{L}_F: \R^{2n+1} \to \R^{2n+1}$ is given by 
$$(\bu, \bv, w)=\mathcal{L}_F(\bx, \bp, z)=(A\bx+\bb, A'\bp+\bb', z+\bc^T\bx+d)$$ 
with $A'=(A^T)^{-1}$, $\bb'=A'\bc$ and $w'=\bv^T\bu-w$, and $\mathcal{L}_F(L_1)=L_2$. 
Then 
\begin{eqnarray*}
&&\mathcal{D}_{L_2} (\mathcal{L}_F(p), \mathcal{L}_F(q))\\
&&\quad =w(p)+w'(q)-\bu(p)^T\bv(q)\\
&&\quad =w(p)-w(q)+\bv(q)^T(\bu(q)-\bu(p))\\
&&\quad =z(p)-z(q)+\bc^T(\bx(p)-\bx(q))+(A'(\bp(q)+\bc))^T(A(\bx(q)-\bx(p)) \qquad \quad\\
&&\quad =z(p)-z(q)+\bp(q)^T(\bx(q)-\bx(p))\\
&&\quad =z(p)+z'(q)-\bx(p)^T\bp(q). \\
&&\quad =\mathcal{D}_{L_1}(p,q). 
\end{eqnarray*}
This completes the proof. 
\qed 

\

Let $(M, \mathcal{U}=\{L_\alpha\})$ be a quasi-Hessian manifold obtained by gluing local models and 
put $\varDelta_M =\{(p, p) \in M \times M\}$. 
Let $U(\varDelta_M)$ denote the subset of $M\times M$ consisting of points $(p, q)$ such that 
there is some local model $L_\alpha$ containing $p, q$. Since $M$ is endowed with the quotient topology, $U(\varDelta_M)$ is an open neighborhood of the diagonal $\varDelta_M$. 

\begin{dfn}\upshape
\label{canonical_divergence}
We set $\mathcal{D}_M(p,q):=\mathcal{D}_{L_\alpha}(p,q)$ at $p, q \in L_\alpha$ for some $\alpha$, then 
$\mathcal{D}_M: U(\varDelta_M) \to \R$ is well-defined by Lemma \ref{invariance}. We call it the {\em canonical divergence} of $M$. 
\end{dfn}

If $M$ is connected and simply connected, then the canonical divergence of $M$ can be extended to the entire space, so we obtain $\mathcal{D}_M:  M\times M \to \R$. 

In Amari-Nagaoka's theory of dually flat structure  \cite{AmariNagaoka00, Amari16}, 
there are two important theorems named by  
{\em extended Pythagorean Theorem} and {\em projection theorem}. 
They take a central role in application to statistical inference, em-algorithm, machine learning and so on. 
These are immediately generalized to our singular setup. 
In the following two theorems, assume that $M$ is a local model (i.e. $M=L \subset \R^{2n+1}$) 
or a connected and simply-connected quasi-Hessian manifold. 
Anyway, the canonical divergence $\mathcal{D}\, (=\mathcal{D}_M)$ is defined on $M \times M$. 

We say that two points $p, q$ are jointed by a curve $c: I \to M$ if there are $t_0, t_1 \in I$ with $c(t_0)=p$ and $c(t_1)=q$. 

\begin{thm} {\bf (Extended Pythagorean Theorem)} \label{pytha_thm}
\label{pythagorean theorem} 
Let $p, q, r \in M$ be three distinct points such that $p$ and $q$ are joined by an $e$-curve $c_e$, and $q$ and $r$ are jointed by an $m$-curve $c_m$, 
and furthermore, $c_e$ and $c_m$ are strictly orthogonal at $q$. Then it holds that 
\begin{align*}
\mathcal{D}(p, r) = \mathcal{D}(p, q) + \mathcal{D}(q, r). 
\end{align*}
\end{thm}

\proof 
Since 
$\mathcal{D}(q,q)=z(q)+z'(q)-\bx(q)^T\bp(q)=0$, we see that 
$$
\mathcal{D}(p, r) - \mathcal{D}(p, q) - \mathcal{D}(q, r)
= -(\bx(p)-\bx(q))^T(\bp(r)-\bp(q)).
$$
The images of the maps $\pi^e_1 \circ c_e$ and $\pi^m_1 \circ c_m$ lie on lines with 
direction vectors, say $\be, \bm$, respectively. 
Then 
$$\bx(p)-\bx(q)=k_0\be, \quad \bp(r)-\bp(q)=k_1\bm$$ 
for some $k_0, k_1 \in \R$. 
The assumption is $\be^T\bm=0$, thus the equality follows. 
\qed

\begin{thm}{\bf (Projection Theorem)} \label{proj_thm}
Let $S$ be a submanifold of $M$ and $c_m : [0, 1] \to L$ an $m$-curve with $q=c_m(1) \in S$. 
Put $p=c_m(0) \in L$. 
Then,  $c_m$ and $S$ are orthogonal at $q$ if and only if 
$q$ is a critical point of the function 
$F=\mathcal{D}(-, p):  S \to \mathbb{R}$. 
The same holds for an $e$-curve $c_e$ and $F=\mathcal{D}(p, -)$. 
\end{thm}

\proof 
Take a generating function $g(\bx_I, \bp_J)$ around $q$. Recall 
$\bp_I= \frac{\rd g}{\rd \bx_I}$, $\bx_J=- \frac{\rd g}{\rd \bp_J}$ and $z=\bp_J^T\bx_J+g(\bx_I, \bp_J)$. 
Let $\gamma=\gamma(s)$ be an immersed curve on $S$ with $\gamma(0)=q$. 
On this curve, 
we have $\frac{d}{ds} g(\bx_I, \bp_J)=\bp_I^T(\frac{d}{ds} \bx_I)-\bx_J^T(\frac{d}{ds} \bp_J)$ and 
$\frac{dz}{ds}=(\frac{d}{ds}\bp_J)^T\bx_J+\bp_J^T(\frac{d}{ds}\bx_J)+\frac{d}{ds}g=\bp^T\frac{d}{ds}\bx$. 
Therefore, we see 
\begin{eqnarray*}
\frac{d(F\circ \gamma)}{ds}(s) 
&=&\frac{d}{ds} (z(\gamma(s))+z'(p)-\bp(p)^T\bx(\gamma(s))) \\
&=& (\bp(\gamma(s))-\bp(p))^T\frac{d(\bx\circ \gamma)}{ds}(s).
\end{eqnarray*}
At $s=0$,  the vector $\bp(q)-\bp(p)$ is a scalar multiple of the direction vector $\bv$ of the line in $\R^n_{\bbp}$ to which the $m$-curve $c_m$ is projected, and $\frac{d}{ds} \bx = d\pi^e_1(\frac{d\gamma}{ds}(0)) \in \R^n_{\bbx}$. Hence, the orthogonality of $S$ and $c_m$ at $q$ is equivalent to that $\frac{d}{ds}F\circ \gamma(0)=0$ for arbitrary $\gamma$, that means that $F$ is critical at $q$. \qed

\begin{exam}\label{Pytha}\upshape
We check the Pythagorean theorem for a toy example in Example \ref{A2}. 
Let 
$$f(\bx)= \frac{x_1^3}{3} + \frac{x_2^2}{2}$$ 
and use  affine local coordinates $\bx=(x_1, x_2)$ for $L$. 
The $m$-Lagrange map is $(x_1, x_2)\mapsto (p_1, p_2)=(x_1^2, x_2)$, and $\Sigma$ is the $x_2$-axis. 
We have already computed the dual potential $z'$, thus for 
$P:=\bx(p)=(a_1, a_2)$ and $Q:=\bx(q)=(b_1, b_2)$, we have 
$$\mathcal{D}(P,Q)=\frac{a_1^3}{3} + \frac{a_2^2}{2}+\frac{2b_1^3}{3}+\frac{b_2^2}{2}-a_1b_1^2-a_2b_2.$$ 
A straight line $p_2=ap_1+b$ on $\R^2_{\bbp}$ corresponds to a parabola $x_2=ax_1^2+b$ on $\R^2_{\bbx}$  (i.e. an $m$-curve). 
Now, for example, take an $m$-curve $c_m$: $x_2=\frac{1}{2}x_1^2$ ($\bm=(2,1)^T$), and two points 
$Q:=(u, \frac{u^2}{2})$ and $R:=(t, \frac{t^2}{2})$ lying on it. 
Take a point $P:=(s, -2(s-u)+\frac{u^2}{2})$ on the straight line on $\R^2_{\bbx}$ (i.e. an $e$-curve) passing through $Q$ directed by $\be=(1,-2)^T$ with $\be^T\bm=0$. 
Then $\triangle PQR$ satisfies the condition, and we see  
$\mathcal{D}(P,Q)+\mathcal{D}(Q, R)=\mathcal{D}(P,R)$. 
It does not matter whether the point $Q$ lies on $\Sigma$ or not. 
\begin{figure}[h]
 \includegraphics[width=9cm, pagebox=cropbox]{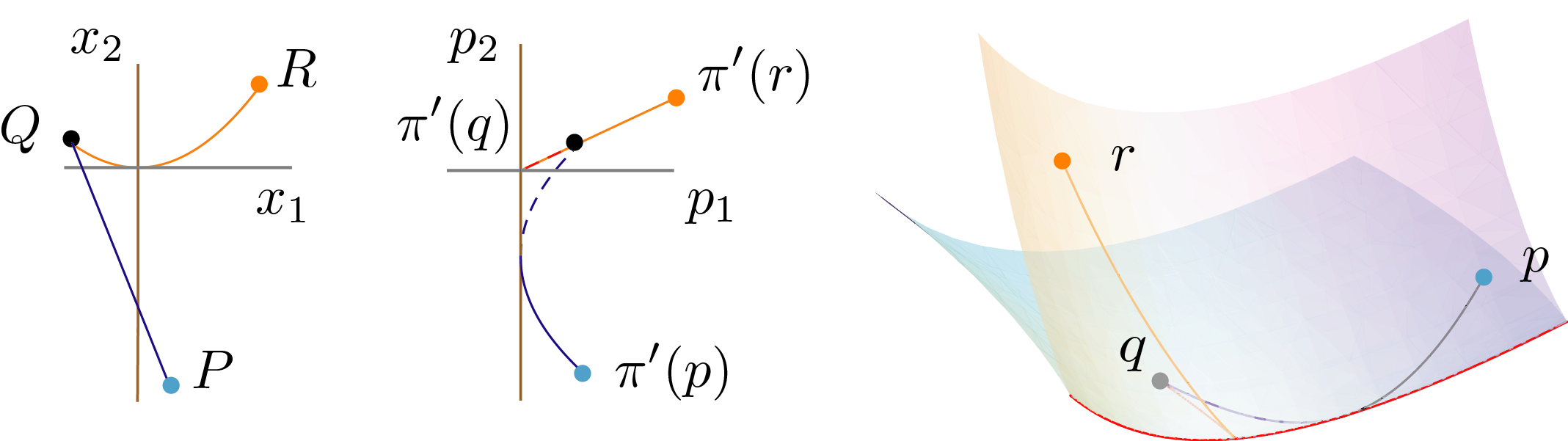}\\
 \caption{\small  
Two projections of the triangle $\triangle pqr$ lying on $L$ of Example \ref{A2}. 
We see a folded triangle on $\R^2_{\bbp}$. 
}\label{fig4}. 
\end{figure}
\end{exam}

\subsection{Weak contrast functions}\label{sec:4-3}
First we recall the definition of {\em contrast functions} (Eguchi \cite{Eguchi}).  
Let $M$ be a manifold and $\rho: U \to \R$ a function defined on 
an open neighborhood $U$ of the diagonal $\mathit\Delta_M \subset M\times M$.  
Given vector fields $X_i\, (1\le i \le k)$, $Y_j \, (1\le j \le l)$ on $M$, 
we set a function 
$$\rho[X_1\cdots X_k|Y_1\cdots Y_l]: M \to \R$$
by assigning to $p \in M$ the value  
$$(X_1)_p\cdots (X_k)_p(Y_1)_q\cdots (Y_l)_q\left(\rho(p, q)\right)|_{p=q}.$$
We also write $\rho[X| - ](r)=X_p\rho(p, q)|_{p=q=r}$ and so on. 
We call $\rho: U \to \R$  a {\em contrast function} of $M$ if it satisfies that 
$$(i)\;\; \rho[-|-]=\rho(p, p)=0 \qquad 
(ii)\;\; \rho[X|-]=\rho[-|X]=0,$$
\begin{center}
(iii) $h(X, Y):=-\rho[X|Y]$ is a pseudo-Riemannian metric on $M$. 
\end{center}
We call $\rho$ a {\em weak contrast function} 
if it satisfies only (i) and (ii). 

Given a contrast function $\rho$,  affine connections are defined by 
$$h(\nabla_XY, Z):=-\rho[XY|Z], \quad h(Y, \nabla^*_XZ):=-\rho[Y|XZ].$$ 
Those connections are torsion-free, mutually dual with respect to $h$, and $\nabla h$ is symmetric, 
and therefore, $(M, h, \nabla)$ becomes a statistical manifold \cite{Eguchi, Matsuzoe}. 
Conversely, given a statistical manifold $M$, one can find a contrast function which reproduces the metric and connections  \cite{Matsumoto} -- it is actually shown in \cite{Matsumoto} 
that for a symmetric $(0,2)$-tensor $h$ (i.e., a possibly degenerate metric) and a symmetric $(0,3)$-tensor $c$, one can find a weak contrast function $\rho: U \to \R$ which satisfies that 
\begin{eqnarray*}
h(X, Y)&=&-\rho[X|Y] \; (=\rho[XY|-]=\rho[-|XY]),  \\
c(X, Y, Z)&=&-\rho[Z|XY]+\rho[XY|Z].
\end{eqnarray*}
Among statistical manifolds, Hessian manifolds admit a notable property: the Bregman divergence is a contrast function, 
and it reproduces the dually flat structure. 
That is extended to our quasi-Hessian manifold and its canonical divergence. 

\begin{thm}\label{contrast}
For a quasi-Hessian manifold $M$, the canonical divergence $\mathcal{D}_M$ is a weak contrast function, and reproduces the quasi-Hessian metric and the canonical cubic tensor by 
\begin{eqnarray*}
h(X, Y)&=&-\mathcal{D}_M[X|Y], \\
 C(X, Y, Z)&=&-\mathcal{D}_M[XY|Z]+\mathcal{D}_M[Z|XY].
\end{eqnarray*}
\end{thm}
\proof 
Since this is a local property, 
take a local model $L_\alpha \subset \R^{2n+1}$. 
Suppose that $g(\bx_I, \bp_J)$ is a generating function 
for $L_\alpha$ around $p \in L_\alpha$. Then $(\bx_I, \bp_J)$ is 
a system of local coordinates for $L_\alpha$ around $p$,  and 
it holds that 
$x_j(q)= - \frac{\rd g}{\rd p_j}(q)$, $p_i(q)=\frac{\rd g}{\rd x_i}(q)$, 
$z(q)=\bp_J(q)^T\bx_J(q)+g(q)$ for 
$q \in L_\alpha$ close to $p$. Hence, 
\begin{eqnarray*}
\mathcal{D}_M(p,q)&=&z(p)-z(q)+\bp(q)^T(\bx(q)-\bx(p))\\
&=&g(p)-g(q)+\bx_J(p)^T(\bp_J(p)-\bp_J(q))+\bp_I(q)^T(\bx_I(q)-\bx_I(p)).
\end{eqnarray*}
Let $\rd_k$ denote $\frac{\rd}{\rd x_k}$ if $k \in I$ and $\frac{\rd}{\rd p_k}$ if $k \in J$, for short. 
Then 
\begin{eqnarray*}
(\rd_k)_p\mathcal{D}_M(p,q)&=&
\epsilon(k)(\rd_kg(p) -\rd_kg(q))+\rd_k\bx_J(p)^T(\bp_J(p)-\bp_J(q)),\\
(\rd_k)_q\mathcal{D}_M(p,q)&=&
(1-\epsilon(k))(\rd_kg(p) -\rd_kg(q))+\rd_k\bp_I(q)^T(\bx_I(q)-\bx_I(p)),
\end{eqnarray*}
where $\epsilon(k)=1$ if $k \in I$ and $0$ if $k \in J$.  
It immediately follows that 
$$\mathcal{D}_M[-|-]=0, \quad  
\mathcal{D}_M[\rd_k|-]=\mathcal{D}_M[-|\rd_k]=0,$$ 
so the divergence is a weak contrast function. 
Put 
$$
\epsilon(k, l)=\left\{
\begin{array}{rl}
1& (k, l \in I)\\
-1& (k, l \in J)\\
0&\mbox{(otherwise)}. 
\end{array}\right.
$$
Then a simple computation shows that 
\begin{eqnarray*}
(\rd_l)_p(\rd_k)_p\mathcal{D}_M(p,q)&=&
\epsilon(k, l)\rd_l\rd_kg(p)+\rd_l\rd_k\bx_J(p)^T(\bp_J(p)-\bp_J(q)),\\
(\rd_l)_q(\rd_k)_q\mathcal{D}_M(p,q)&=&
\epsilon(k, l)\rd_l\rd_kg(q)+\rd_l\rd_k\bp_I(q)^T(\bx_I(q)-\bx_I(p)). \qquad\quad
\end{eqnarray*}
Hence 
$\mathcal{D}_M[\rd_k\rd_l|-]= h(\rd_k, \rd_l)$ 
and 
$$\mathcal{D}_M[\rd_k\rd_l|\rd_m]-\mathcal{D}_M[\rd_m|\rd_k\rd_l]=- \rd_k\rd_l\rd_mg$$ 
for any $k, l, m$. This coincides with the cubic tenser $C$ by Proposition \ref{C_local} up to the sign. 
\qed

\section{Discussions}\label{sec:5}

We shortly discuss possible directions or proposals for further researches. 

\subsection{Pre-Frobenius structure}\label{sec:5-1} 
In mathematical physics such as string theory, there often arise manifolds endowed with commutative and associative multiplication on tangent spaces satisfying certain properties, called (several variations of)  {\em Frobenius manifolds} \cite{Dubrovin}. 
Now, let $(M, h, C)$ be a {\em flat} Hessian manifold, i.e., the metric connection with respect to $h$ is flat. 
Then $M$ naturally carries a (weak) version of Frobenius structure  \cite[\S 2]{Totaro}. 
Put $C_{ijk}=C(\rd_i, \rd_j, \rd_k)$ using $\nabla$-affine coordinates, 
and we may take them as structure constants 
to define a multiplication on $T_pM$: 
$$\rd_i \circ \rd_j := \sum_{k,l} C_{ijk}h^{kl} \rd_l.$$
Since $C$ is symmetric, it is commutative. 
The associativity, $(\rd_i \circ \rd_j)\circ \rd_k=\rd_i \circ (\rd_j\circ \rd_k)$, 
is written down to 
$$\sum_{a,b}(C_{ijb}C_{kla}-C_{ila}C_{jkb})=0 \qquad (\forall\, i, j, k, l),$$
and a bit surprisingly,  the left hand side coincides with the curvature tensor for the Levi-Civita connection of $h$  \cite{Duistermaat, Kito}; 
the equation is actually known as the WDVV equation in string theory. 
Moreover, it is easy to see that 
the multiplication is compatible with the metric: $h(\rd_i \circ \rd_j, \rd_k)=h(\rd_i, \rd_j \circ \rd_k)$. 
Then the tuple $(M, h, \circ)$ becomes a weak pre-Frobenius manifold (cf. \cite{Dubrovin, Hertling}). 
For a quasi-Hessian manifold $M$, the symmetric cubic tensor $C=(C_{ijk})$ is defined everywhere, but $h^{kl}$ is not; 
even though, the WDVV equation makes sense.  Then, at least for every $p \in \Sigma$ (pointwise),  
the quotient $T_pM/{\rm null}(h_p)$ carries a Frobenius algebra structure. 

A new pre-Frobenius structure on  a certain space of probability distributions has recently been found using the Hessian geometry on convex cones and paracomplex structure in \cite{Combe}. 
Also from the context of Poisson and paraK\"ahler geometry, the notion of contravariant pseudo-Hessian manifolds has been introduced in \cite{Boucetta19}, which is actually very close to our quasi-Hessian manifolds 
with degenerate potentials.  Those should be mutually related. 

As a different question from the above, 
more interesting is local geometry of quasi-Hessian $M$ in relation with the Saito-Givental theory  --  under a certain condition, the germ of $M$ at a point should be a real geometry counterpart to analytic spectrum of a massive F-manifold 
(cf. \cite[\S 3]{Hertling}; the analytic spectrum is a certain holomorphic Legendre submanifold of $\C^{2n+1}$ defined by a versal deformation of a complex isolated hypersurface singularity as its a generating family). Perhaps, this was essentially posed by Arnol'd \cite[\S4]{Arnold90}.

\subsection{Statistical inference and machine learning}\label{sec:5-2}
\begin{figure}[h]
 \includegraphics[width=3cm, pagebox=cropbox]{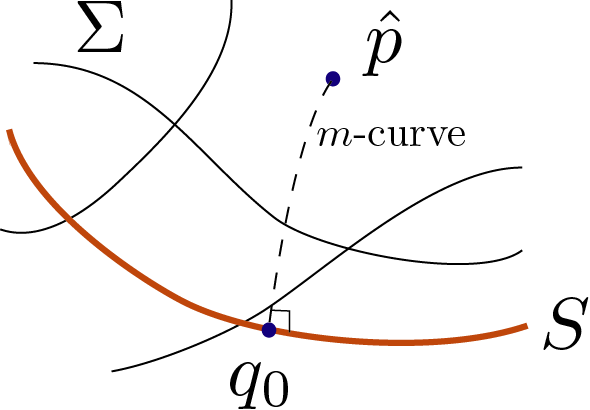}\\
 \caption{\small A conceptual figure for statistical inference.}\label{fig5}
\end{figure}

Suppose that our statistical model $S$ is a curved exponential family, i.e., a submanifold of an exponential family $M$ (see  Example \ref{exponential_family}). 
Let $\mathcal{D}: M\times M \to \R$ be the associated Bregman divergence, which is known to coincide 
with the Kullback-Leibler divergence $$\mathcal{D}_{KL}(q, p)=\int q(u)\log\frac{q(u)}{p(u)} du$$ 
measuring an `asymmetrical distance' from a distribution (density function) $q=q(u)$ to another $p=p(u)$. A given data set $\{{\bf u}_i\}$ yields an observed point $\hat{p} \in M$, then the task of statistical inference is to find $q_0 \in S$ which best approximates the point $\hat{p}$. 
Information geometry \cite{AmariNagaoka00, Amari16} provides a clear geometric understanding on the maximum likelihood estimate (MLE), 
that is,  the MLE assigns to $\hat{p}$ the point $q_0\in S$ which attains the minimum of $\mathcal{D}(\cdot, \hat{p}): S \to \R$, and especially, $\hat{p}$ is projected to $q_0$ along an $m$-geodesic ($m$-curve) being orthogonal to $S$ at $q_0$.  
We have shown that this assertion is valid even in case that $M$ admits the locus $\Sigma$ where the Fisher-Rao metric is degenerate (Theorem \ref{proj_thm}), see Fig.\,\ref{fig5}. 

If $\hat{p}$ is sufficiently close to $S$ and far from $\Sigma$, then the asymptotic theory of estimation is discussed. However, in practice, we may not be able to know if $\hat{p}$ is the case. 
For instance, it often happens that the MLE has multiple local minimums, i.e., the maximum likelihood equation may have multiple roots. Then, as $\hat{p}$ varies by renewing the data,  {\em catastrophe phenomena}  -- the birth and death of min/max. points --  can happen. Actually, the ambiguity of root selection in MLE has been studied in practical and numerical approach (cf.   \cite[\S 4]{Small}), while there seems to be less theoretical approach so far. 
Our framework provides a right way from information geometry. Define 
$$F: S \times M \to \R\qquad F(q, p):=\mathcal{D}(\iota(q), p)$$
and we may consider $F$ as a global generating family \cite[p.323]{AGV}, i.e., it defines a Legendre submanifold of $T^*M\times \R$ by 
$$L_S:=\left\{\; (p, \eta, z)  \; \middle|  \; \exists\, q \in S, \; \frac{\rd F}{\rd q}(q, p)=0, \; \eta=\frac{\rd F}{\rd p}(q, p), \; z=F(p, \eta) \right\}$$
where we roughly denote by $\frac{\rd F}{\rd q}$ the differential with respect to $S$ and so on. 
This gives a typical example of a quasi-Hessian manifold. 
The critical value set of the Lagrange map $\pi: L_S \to M$ ($\pi(p,\eta, z) \mapsto p$) is nothing but the envelope of the family of all $m$-curves on $M$ which are orthogonal to $S$;  we call it the {\em $m$-caustics determined by $S$}. 
If $S$ is not $\nabla$-flat,  the $m$-caustics usually appear (that reflects the $\nabla^*$-extrinsic geometry of $S$ in $M$). 
It turns out that the catastrophe phenomenon mentioned above arises when the data manifold $D$ intersects with the $m$-caustics determined by $S$. Conversely, for a given data manifold $D$, we may consider the restriction of $\mathcal{D}$ to $M\times D$ and define the {\em $e$-caustics determined by $D$} similarly. Interaction between these two $e/m$-caustics can be involved and affect the performance of EM-algorithm (cf. Amari \cite[Chap. 8]{Amari16}). 
Note that in principle, the above strategy may be adapted to any divergence and any statistical model. 
The detail will be discussed somewhere else.

As described in Amari \cite[Chap.11]{Amari16}, a class of learning machines is also based on the Bregman divergence $\mathcal{D}_\phi$ of convex functions $\phi$. Now, as an attempt, suppose that $\phi$ is a nonconvex function (possibly with inflection points).  Read $\mathcal{D}_\phi$ to be the corresponding canonical divergence in our sense (see \S \ref{sec:4-2}).  
Here we would like to notice that the same proofs in convex case do often work to obtain slightly weaker results for such general $\phi$ --  an easy example is Theorem 11.1 of \cite{Amari16}, which is read off as  ``{\em the $k$-mean $\eta_C:=\frac{1}{k}\sum x_i$ of a cluster $C=\{x_i\}_{i=1}^k$ in $\R^n$ is always a critical point of $\mathcal{D}_\phi(C, -):=\frac{1}{k}\sum\mathcal{D}_\phi(x_i, -)$, and all other critical points are obtained from $\eta_C$ and $\ker \nabla^2\phi$}". 
We expect a similar result for some other optimization algorithm. 
On the other hand, almost all statistical learning machines allow Fisher-Rao matrices to be degenerate \cite{FK, Watanabe}. In particular, as in \cite[Chap.12]{Amari16}, most of deep learning machines use the Gaussian noise with a fixed (co)variance for regression; then the parameter space $M$ becomes a self-dual Riemannian manifold $(h, \nabla=\nabla^*)$ off the degeneracy locus $\Sigma$ of $h$ having many components.  
We seek another scheme for measuring errors which is compatible with our singular model.

\subsection{Conclusion}\label{sec:5-3}
In the present paper, we have proposed an information geometry for singular models 
from the viewpoint of contact geometry and singularity theory. 
We have introduced quasi-Hessian manifolds, 
which extend the notion of dually flat manifolds of Amari-Nagaoka so that the Hessian metric can be degenerate, but the canonical cubic tensor is consistently defined on the entire space. 
Most notable is that the extended Pythagorean theorem and projection theorem are valid even in this singular setup. 

There are several further directions as mentioned above. 
We end by adding a few more comments.  
There is an on-going project of the first author on local classification of singularities of $em$-wavefronts in flat affine coordinates, which extends an old work of Ekeland \cite{Ekeland} in nonconvex optimization and leads to affine differential geometry of wavefronts (cf. \cite{SUY}). 
Secondly, since a quasi-Hessian manifold is embedded in some contact manifold, 
we may think of the Hamiltonian-Jacobi method for time evolution of quasi-Hessian manifolds (wavefront propagation) and semi-classical quantization (WKB analysis) in our framework (cf. \cite{Arnold89}). 
Finally, it would be valuable to find some connections with preceding excellent works on singular statistical models \cite{Amari16, FK, Watanabe} -- especially, we hope that the theory of singular Legendre varieties and Legendre currents would make a bridge between the differential geometric method \cite{AmariNagaoka00, Amari16} and the  algebro-geometric method \cite{Watanabe}.


\end{document}